\documentclass{amsart}

\usepackage{amssymb}
\usepackage{amsmath}
\usepackage{amsthm}
\usepackage{amsfonts}

\usepackage{a4wide}

\newtheorem{theorem}{Theorem}

\theoremstyle{definition}

\begin{document}


\title[Periodic Schur \(\sigma\)-groups]
{Periodic Schur \(\sigma\)-groups of \\ non-elementary bicyclic type}

\author{Daniel C. Mayer}
\address{Naglergasse 53\\8010 Graz\\Austria}
\email{algebraic.number.theory@algebra.at}
\urladdr{http://www.algebra.at}

\thanks{Research supported by the Austrian Science Fund (FWF): projects J0497-PHY, P26008-N25, and by EUREA}



\date{Tuesday, 26 October 2021 (Austrian National Celebration Day)}


\begin{abstract}
Infinitely many large Schur \(\sigma\)-groups \(G\) with
non-elementary bicyclic commutator quotient \(G/G^\prime\simeq C_{3^e}\times C_3\), \(e\ge 2\),
are constructed as periodic sequences of vertices in descendant trees of finite \(3\)-groups.
A single root gives rise to pairs of metabelian groups \(G\)
with logarithmic order \(\mathrm{lo}(G)=4+e\) for \(e\ge 3\).
Three roots are ancestors of pairs of non-metabelian groups \(G\) with
moderate rank distribution \(\varrho(G)\sim (2,2,3;3)\) and
\(\mathrm{lo}(G)=7+e\) for \(e\ge 5\).
Twentyseven roots produce sextets of non-metabelian groups \(G\) with
elevated rank distribution \(\varrho(G)=(3,3,3;3)\) and
\(\mathrm{lo}(G)=19+e\) for \(e\ge 9\).
The soluble length of non-metabelian groups is always \(\mathrm{sl}(G)=3\).
The groups can be realized as \(3\)-class field tower groups \(\mathrm{Gal}(\mathrm{F}_3^\infty(K)/K)\)
of imaginary quadratic number fields \(K=\mathbb{Q}(\sqrt{d})\)
with fundamental discriminants \(d<0\).
\end{abstract}

\maketitle


\section{Introduction}
\label{s:Intro}

\noindent
Why are Schur \(\sigma\)-groups of eminent algebraic and arithmetic relevance?
For an assigned odd prime number \(p\ge 3\),
the automorphism group \(G_\infty=\mathrm{Gal}(\mathrm{F}_p^\infty(K)/K)\)
of the maximal unramified pro-\(p\) extension \(\mathrm{F}_p^\infty(K)\)
of an imaginary quadratic number field \(K=\mathbb{Q}(\sqrt{d})\)
with fundamental discriminant \(d<0\) is a Schur \(\sigma\)-group
\cite{Sh1964,KoVe1975,Ag1998,BuMa2015,BBH2017},
in fact, either a finite group with order a power of \(p\)
or an infinite pro-\(p\) group.
All Galois groups \(G_n=\mathrm{Gal}(\mathrm{F}_p^n(K)/K)\)
of the stages of the unramified Hilbert \(p\)-class field tower
\[
K\le\mathrm{F}_p^1(K)\le\mathrm{F}_p^2(K)\le\ldots\le\mathrm{F}_p^n(K)\le\ldots\le\mathrm{F}_p^\infty(K)
\]
of \(K\) are derived quotients \(G_n\simeq G_\infty/G_\infty^{(n)}\)
of the Schur \(\sigma\)-group \(G_\infty\).
In turn, imaginary quadratic fields \(K\) are the simplest algebraic number fields
with smallest possible degree \(2\) and trivial torsion-free unit group \(U_K/W_K\).
At least for the smallest odd prime number \(p=3\),
imaginary quadratic fields are distinguished by the experimental constructibility
of their unramified abelian extensions
\cite{Fi2001,BCP1997,BCFS2021,MAGMA2021}
with relative degrees \(p\) and usually also \(p^2\)
in a few minutes of CPU-time.

In our most recent investigations,
we succeeded in finding unexpected periodic sequences
of Schur \(\sigma\)-groups \(G\) possessing a bicyclic commutator quotient
\(G/G^\prime\simeq C_{3^e}\times C_3\)
with one non-elementary component and logarithmic exponent \(e\ge 2\).
Periodicity sets in for a minimal exponent \(e\ge e_0\)
in dependence on the particular kind of the Schur \(\sigma\)-groups \(G\).
The value \(e_0\) is given by the \(p\)-nilpotency class \(\mathrm{cl}_p(G)\) of \(G\).
There exist periodic sequences of
\begin{itemize}
\item
metabelian Schur \(\sigma\)-groups \(G\)
with logarithmic order \(\mathrm{lo}(G)=4+e\) for \(e\ge e_0=3\)
\cite{Ma2021a},
\item
non-metabelian Schur \(\sigma\)-groups \(G\)
with moderate rank distribution \(\varrho(G)\sim (2,2,3;3)\) or \(\varrho(G)\sim (2,2,2;3)\)
and logarithmic order \(\mathrm{lo}(G)=7+e\) for \(e\ge e_0=5\)
\cite{Ma2021a,Ma2021c},
\item
non-metabelian Schur \(\sigma\)-groups \(G\)
with elevated rank distribution \(\varrho(G)=(3,3,3;3)\)
and logarithmic order \(\mathrm{lo}(G)=19+e\) for \(e\ge e_0=9\)
\cite{Ma2021b}.
\end{itemize}


\section{Metabelian Schur-groups \(G\) with \(G/G^\prime\simeq (3^e,3)\), \(e\ge 1\)}
\label{s:Metabelian1}

\noindent
A common feature of all periodic constructions in this article
is the constitution of the resulting \(p\)-descendant tree by
\begin{itemize}
\item
an \textit{infinite main trunk} of vertices \(T\) with parametrized presentation,
\item
graph theoretically isomorphic \textit{finite twigs} emanating from each vertex of the main trunk.
\end{itemize}
The leaves \(G\) of the twigs are usually Schur \(\sigma\)-groups.
In order to avoid ambiguity, we use subscripts \(T_e,G_e\) in all the following theorems.

Since the simplest main trunk consists of all abelian groups of type \((3^e,3)\),
we begin with this important instance,
although the leaves of the twigs are only Schur-groups (with balanced presentation)
but not \(\sigma\)-groups (without generator-inverting automorphism).

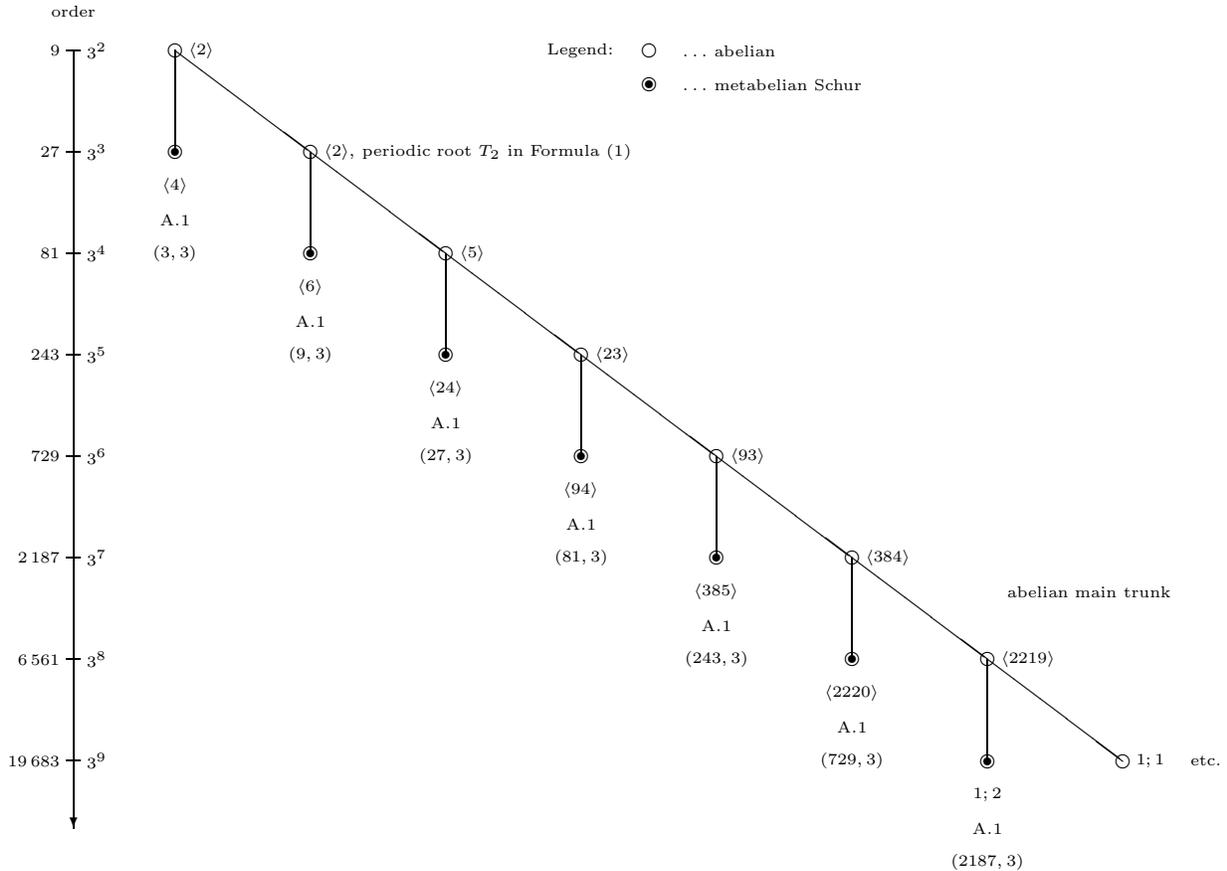
\begin{figure}[hb]
\caption{Periodic metabelian Schur-groups \(G\) with \(G/G^\prime\simeq (3^e,3)\), \(e\ge 2\)}
\label{fig:Schur2}

{\tiny

\setlength{\unitlength}{0.9cm}
\begin{picture}(16,12.5)(-10,-6.5)

\put(-10,5.5){\makebox(0,0)[cb]{order}}

\put(-10,5){\line(0,-1){10.5}}
\multiput(-10.1,5)(0,-1.5){8}{\line(1,0){0.2}}

\put(-10.2,5){\makebox(0,0)[rc]{\(9\)}}
\put(-9.8,5){\makebox(0,0)[lc]{\(3^2\)}}
\put(-10.2,3.5){\makebox(0,0)[rc]{\(27\)}}
\put(-9.8,3.5){\makebox(0,0)[lc]{\(3^3\)}}
\put(-10.2,2){\makebox(0,0)[rc]{\(81\)}}
\put(-9.8,2){\makebox(0,0)[lc]{\(3^4\)}}
\put(-10.2,0.5){\makebox(0,0)[rc]{\(243\)}}
\put(-9.8,0.5){\makebox(0,0)[lc]{\(3^5\)}}
\put(-10.2,-1){\makebox(0,0)[rc]{\(729\)}}
\put(-9.8,-1){\makebox(0,0)[lc]{\(3^6\)}}
\put(-10.2,-2.5){\makebox(0,0)[rc]{\(2\,187\)}}
\put(-9.8,-2.5){\makebox(0,0)[lc]{\(3^7\)}}
\put(-10.2,-4){\makebox(0,0)[rc]{\(6\,561\)}}
\put(-9.8,-4){\makebox(0,0)[lc]{\(3^8\)}}
\put(-10.2,-5.5){\makebox(0,0)[rc]{\(19\,683\)}}
\put(-9.8,-5.5){\makebox(0,0)[lc]{\(3^9\)}}

\put(-10,-5.5){\vector(0,-1){1}}

\put(-8.5,5){\circle{0.2}}

\put(-8.5,3.5){\circle{0.2}}
\put(-8.5,3.5){\circle*{0.1}}

\put(-6.5,3.5){\circle{0.2}}
\put(-6.5,2){\circle{0.2}}
\put(-6.5,2){\circle*{0.1}}

\put(-4.5,2){\circle{0.2}}
\put(-4.5,0.5){\circle{0.2}}
\put(-4.5,0.5){\circle*{0.1}}

\put(-2.5,0.5){\circle{0.2}}
\put(-2.5,-1){\circle{0.2}}
\put(-2.5,-1){\circle*{0.1}}

\put(-0.5,-1){\circle{0.2}}
\put(-0.5,-2.5){\circle{0.2}}
\put(-0.5,-2.5){\circle*{0.1}}

\put(1.5,-2.5){\circle{0.2}}
\put(1.5,-4){\circle{0.2}}
\put(1.5,-4){\circle*{0.1}}

\put(3.5,-4){\circle{0.2}}
\put(3.5,-5.5){\circle{0.2}}
\put(3.5,-5.5){\circle*{0.1}}

\put(5.5,-5.5){\circle{0.2}}

\put(-8.5,5){\line(0,-1){1.5}}
\put(-8.5,5){\line(4,-3){2}}

\put(-6.5,3.5){\line(0,-1){1.5}}
\put(-6.5,3.5){\line(4,-3){2}}

\put(-4.5,2){\line(0,-1){1.5}}
\put(-4.5,2){\line(4,-3){2}}

\put(-2.5,0.5){\line(0,-1){1.5}}
\put(-2.5,0.5){\line(4,-3){2}}

\put(-0.5,-1){\line(0,-1){1.5}}
\put(-0.5,-1){\line(4,-3){2}}

\put(1.5,-2.5){\line(0,-1){1.5}}
\put(1.5,-2.5){\line(4,-3){2}}

\put(3.5,-4){\line(0,-1){1.5}}
\put(3.5,-4){\line(4,-3){2}}


\put(-8.3,5){\makebox(0,0)[lc]{\(\langle 2\rangle\)}}
\put(-8.5,3){\makebox(0,0)[cc]{\(\langle 4\rangle\)}}
\put(-8.5,2.5){\makebox(0,0)[cc]{\(\mathrm{A}.1\)}}
\put(-8.5,2){\makebox(0,0)[cc]{\((3,3)\)}}

\put(-6.3,3.5){\makebox(0,0)[lc]{\(\langle 2\rangle\), periodic root \(T_2\) in Formula \eqref{eqn:Abelian}}}
\put(-6.5,1.5){\makebox(0,0)[cc]{\(\langle 6\rangle\)}}
\put(-6.5,1){\makebox(0,0)[cc]{\(\mathrm{A}.1\)}}
\put(-6.5,0.5){\makebox(0,0)[cc]{\((9,3)\)}}

\put(-4.3,2){\makebox(0,0)[lc]{\(\langle 5\rangle\)}}
\put(-4.5,0){\makebox(0,0)[cc]{\(\langle 24\rangle\)}}
\put(-4.5,-0.5){\makebox(0,0)[cc]{\(\mathrm{A}.1\)}}
\put(-4.5,-1){\makebox(0,0)[cc]{\((27,3)\)}}

\put(-2.3,0.5){\makebox(0,0)[lc]{\(\langle 23\rangle\)}}
\put(-2.5,-1.5){\makebox(0,0)[cc]{\(\langle 94\rangle\)}}
\put(-2.5,-2){\makebox(0,0)[cc]{\(\mathrm{A}.1\)}}
\put(-2.5,-2.5){\makebox(0,0)[cc]{\((81,3)\)}}

\put(-0.3,-1){\makebox(0,0)[lc]{\(\langle 93\rangle\)}}
\put(-0.5,-3){\makebox(0,0)[cc]{\(\langle 385\rangle\)}}
\put(-0.5,-3.5){\makebox(0,0)[cc]{\(\mathrm{A}.1\)}}
\put(-0.5,-4){\makebox(0,0)[cc]{\((243,3)\)}}

\put(1.7,-2.5){\makebox(0,0)[lc]{\(\langle 384\rangle\)}}
\put(1.5,-4.5){\makebox(0,0)[cc]{\(\langle 2220\rangle\)}}
\put(1.5,-5){\makebox(0,0)[cc]{\(\mathrm{A}.1\)}}
\put(1.5,-5.5){\makebox(0,0)[cc]{\((729,3)\)}}

\put(3.7,-4){\makebox(0,0)[lc]{\(\langle 2219\rangle\)}}
\put(3.5,-6){\makebox(0,0)[cc]{\(1;2\)}}
\put(3.5,-6.5){\makebox(0,0)[cc]{\(\mathrm{A}.1\)}}
\put(3.5,-7){\makebox(0,0)[cc]{\((2187,3)\)}}

\put(5.7,-5.5){\makebox(0,0)[lc]{\(1;1\)}}

\put(5,-3){\makebox(0,0)[cc]{abelian main trunk}}
\put(6.5,-5.5){\makebox(0,0)[lc]{etc.}}


\put(-3,5){\makebox(0,0)[lc]{Legend:}}

\put(-1.5,5){\circle{0.2}}
\put(-1,5){\makebox(0,0)[lc]{\(\ldots\) abelian}}

\put(-1.5,4.5){\circle{0.2}}
\put(-1.5,4.5){\circle*{0.1}}
\put(-1,4.5){\makebox(0,0)[lc]{\(\ldots\) metabelian Schur}}

\end{picture}

}

\end{figure}

\begin{theorem}
\label{thm:AbelianChain}
For each logarithmic exponent \(e\ge 2\),
the \textbf{unique} non-elementary bi-heterocyclic \(3\)-group \(T_e\) of type \((3^e,3)\hat{=}(e1)\)
is given by the \textbf{periodic sequence} of iterated \(p\)-descendants
\begin{equation}
\label{eqn:Abelian}
T_e\simeq\mathrm{SmallGroup}(27,2)(-\#1;1)^{e-2}.
\end{equation}
These abelian groups have logarithmic order \(\mathrm{lo}(T_e)=1+e\),
punctured transfer kernel type \(\mathrm{a}.1\), \(\varkappa(T_e)=(000;0)\), and
first abelian quotient invariants \(\alpha_1(T_e)\sim (e,e,e;(e-1)1)\).
They form the infinite main trunk of a descendant tree
with singlets as finite twigs.
For each
\(e\ge 2\),
the singlet
\begin{equation}
\label{eqn:A1}
G_e\simeq\mathrm{SmallGroup}(27,2)(-\#1;1)^{e-2}-\#1;2
\end{equation}
is the \textbf{unique} metabelian Schur-group \(G_e\) with
commutator quotient \(G_e/G_e^\prime\simeq (3^e,3)\hat{=}(e1)\) and
punctured transfer kernel type \(\mathrm{A}.1\), \(\varkappa(G_e)\sim (111;1)\).
It has first abelian quotient invariants \(\alpha_1(G_e)\sim (e+1,e+1,e+1;e1)\),
\(\mathrm{lo}(G_e)=2+e\),
and is not a \(\sigma\)-group.
See Figure
\ref{fig:Schur2}.
\end{theorem}


\begin{proof}
Let \(s_2=\lbrack y,x\rbrack\) denote the main commutator and
\(w=x^{3^e}\) the last non-trivial power.
For each \(e\ge 1\), we have parametrized presentations
\begin{equation}
\label{eqn:Pres1}
\begin{aligned}
T_{e+1} &= \langle x,y,w\mid x^{3^{e+1}}=w^3=1,\ y^3=1\rangle, \\
G_e &= \langle x,y,s_2,w\mid x^{3^{e+1}}=w^3=1,\ y^3=1,\ s_2=w\rangle.
\end{aligned}
\end{equation}
For \(e\ge 1\), the last non-trivial lower \(p\)-central is given by
\(P_{e}(G_e)=\langle w\rangle\) respectively \(P_{e}(T_{e+1})=\langle w\rangle\),
whence both groups share the common \(p\)-parent
\(G_e/P_{e}(G_e)\simeq T_{e+1}/P_{e}(T_{e+1})\simeq T_e\).
Observe that the \(p\)-nilpotency class is given by
\(\mathrm{cl}_p(G_e)=e+1\) respectively \(\mathrm{cl}_p(T_{e+1})=e+1\).

The \(p\)-group generation algorithm
\cite{HEO2005}
by Newman
\cite{Nm1977}
and O'Brien
\cite{Ob1990}
is implemented in the ANUPQ package
\cite{GNO2006}
of the computational algebra system Magma
\cite{MAGMA2021,BCFS2021,BCP1997}.
This algorithm is used to construct
all immediate \(p\)-descendants of an assigned finite \(p\)-group.
Vertices are identified by
absolute counters \(\mathrm{SmallGroup}(o,i)=\langle o,i\rangle=\langle i\rangle\),
defined in the SmallGroups database
\cite{BEO2005}
for orders \(o\le 3^8\),
or by relative counters \(P-\#s;i\) with respect to the \(p\)-parent \(P\) and the step size \(s\),
assigned by the ANUPQ package
\cite{GNO2006}
for \(o\ge 3^9\).
Repeated recursive applications of the algorithm
eventually produce Figure
\ref{fig:Schur2},
and thus confirm Formula
\eqref{eqn:Abelian}
and
\eqref{eqn:A1}.

Special care is required for the tree root \(T_1=\langle 9,2\rangle\) only.
It has exceptional nuclear rank \(n(T_1)=3\),
and we can neglect step sizes \(s\in\lbrace 2,3\rbrace\).
Among the immediate \(p\)-descendants of step size \(s=1\), the abelian group
\(T_2=\langle 27,2\rangle\) is exo-genetic with regular identifier \(T_1-\#1;1\), but
\(G_1=\langle 27,4\rangle\) is endo-genetic with irregular identifier \(T_1-\#1;3\),
since \(T_1-\#1;2\simeq\langle 27,3\rangle\)
is the extra-special \(3\)-group with punctured transfer kernel type \(\mathrm{a}.1\), \(\varkappa=(000;0)\).
Therefore, it is not possible to extend Theorem
\ref{thm:AbelianChain}
to all \(e\ge 1\) in a uniform way.
\end{proof}


\section{Metabelian Schur \(\sigma\)-groups \(G\) with \(G/G^\prime\simeq (3^e,3)\), \(e\ge 3\)}
\label{s:Metabelian3}

\begin{figure}[hb]
\caption{Periodic metabelian Schur \(\sigma\)-groups \(G\) with \(G/G^\prime\simeq (3^e,3)\), \(e\ge 3\)}
\label{fig:SchurSigma3}

{\tiny

\setlength{\unitlength}{0.9cm}
\begin{picture}(16,11)(-10,-8)

\put(-10,2.5){\makebox(0,0)[cb]{order}}

\put(-10,2){\line(0,-1){9}}
\multiput(-10.1,2)(0,-1.5){7}{\line(1,0){0.2}}

\put(-10.2,2){\makebox(0,0)[rc]{\(243\)}}
\put(-9.8,2){\makebox(0,0)[lc]{\(3^5\)}}
\put(-10.2,0.5){\makebox(0,0)[rc]{\(729\)}}
\put(-9.8,0.5){\makebox(0,0)[lc]{\(3^6\)}}
\put(-10.2,-1){\makebox(0,0)[rc]{\(2\,187\)}}
\put(-9.8,-1){\makebox(0,0)[lc]{\(3^7\)}}
\put(-10.2,-2.5){\makebox(0,0)[rc]{\(6\,561\)}}
\put(-9.8,-2.5){\makebox(0,0)[lc]{\(3^8\)}}
\put(-10.2,-4){\makebox(0,0)[rc]{\(19\,683\)}}
\put(-9.8,-4){\makebox(0,0)[lc]{\(3^9\)}}
\put(-10.2,-5.5){\makebox(0,0)[rc]{\(59\,049\)}}
\put(-9.8,-5.5){\makebox(0,0)[lc]{\(3^{10}\)}}
\put(-10.2,-7){\makebox(0,0)[rc]{\(177\,147\)}}
\put(-9.8,-7){\makebox(0,0)[lc]{\(3^{11}\)}}

\put(-10,-7){\vector(0,-1){1}}

\put(-9,2){\circle{0.2}}
\put(-9,2){\circle*{0.1}}
\put(-8,2){\circle{0.2}}
\put(-8,2){\circle*{0.1}}

\put(-7,0.5){\circle{0.2}}
\put(-7,0.5){\circle*{0.1}}
\put(-6,0.5){\circle{0.2}}
\put(-6,0.5){\circle*{0.1}}

\put(-4.5,0.5){\circle{0.2}}
\put(-5,-1){\circle{0.2}}
\put(-5,-1){\circle*{0.1}}
\put(-4,-1){\circle{0.2}}
\put(-4,-1){\circle*{0.1}}

\put(-2.5,-1){\circle{0.2}}
\put(-3,-2.5){\circle{0.2}}
\put(-3,-2.5){\circle*{0.1}}
\put(-2,-2.5){\circle{0.2}}
\put(-2,-2.5){\circle*{0.1}}

\put(-0.5,-2.5){\circle{0.2}}
\put(-1,-4){\circle{0.2}}
\put(-1,-4){\circle*{0.1}}
\put(0,-4){\circle{0.2}}
\put(0,-4){\circle*{0.1}}

\put(1.5,-4){\circle{0.2}}
\put(1,-5.5){\circle{0.2}}
\put(1,-5.5){\circle*{0.1}}
\put(2,-5.5){\circle{0.2}}
\put(2,-5.5){\circle*{0.1}}

\put(3.5,-5.5){\circle{0.2}}
\put(3,-7){\circle{0.2}}
\put(3,-7){\circle*{0.1}}
\put(4,-7){\circle{0.2}}
\put(4,-7){\circle*{0.1}}


\put(-4.5,0.5){\line(-1,-3){0.5}}
\put(-4.5,0.5){\line(1,-3){0.5}}
\put(-4.5,0.5){\line(4,-3){2}}

\put(-2.5,-1){\line(-1,-3){0.5}}
\put(-2.5,-1){\line(1,-3){0.5}}
\put(-2.5,-1){\line(4,-3){2}}

\put(-0.5,-2.5){\line(-1,-3){0.5}}
\put(-0.5,-2.5){\line(1,-3){0.5}}
\put(-0.5,-2.5){\line(4,-3){2}}

\put(1.5,-4){\line(-1,-3){0.5}}
\put(1.5,-4){\line(1,-3){0.5}}
\put(1.5,-4){\line(4,-3){2}}

\put(3.5,-5.5){\line(-1,-3){0.5}}
\put(3.5,-5.5){\line(1,-3){0.5}}

\put(3.5,-5.5){\line(4,-3){2}}


\put(-9,1.5){\makebox(0,0)[cc]{\(\langle 5\rangle\)}}
\put(-8,1.5){\makebox(0,0)[cc]{\(\langle 7\rangle\)}}
\put(-9,1){\makebox(0,0)[cc]{\(\mathrm{D}.10\)}}
\put(-8,1){\makebox(0,0)[cc]{\(\mathrm{D}.5\)}}
\put(-8.5,0.5){\makebox(0,0)[cc]{\((3,3)\)}}

\put(-7,0){\makebox(0,0)[cc]{\(\langle 14\rangle\)}}
\put(-6,0){\makebox(0,0)[cc]{\(\langle 15\rangle\)}}
\put(-6.5,-0.5){\makebox(0,0)[cc]{\(\mathrm{D}.11\)}}
\put(-6.5,-1){\makebox(0,0)[cc]{\((9,3)\)}}

\put(-4.3,0.5){\makebox(0,0)[lc]{\(\langle 8\rangle\), periodic root \(T_3\) in Formula \eqref{eqn:b16}}}
\put(-5,-1.5){\makebox(0,0)[cc]{\(\langle 121\rangle\)}}
\put(-4,-1.5){\makebox(0,0)[cc]{\(\langle 122\rangle\)}}
\put(-4.5,-2){\makebox(0,0)[cc]{\(\mathrm{D}.11\)}}
\put(-4.5,-2.5){\makebox(0,0)[cc]{\((27,3)\)}}

\put(-2.3,-1){\makebox(0,0)[lc]{\(\langle 120\rangle\)}}
\put(-3,-3){\makebox(0,0)[cc]{\(\langle 975\rangle\)}}
\put(-2,-3){\makebox(0,0)[cc]{\(\langle 976\rangle\)}}
\put(-2.5,-3.5){\makebox(0,0)[cc]{\(\mathrm{D}.11\)}}
\put(-2.5,-4){\makebox(0,0)[cc]{\((81,3)\)}}

\put(-0.3,-2.5){\makebox(0,0)[lc]{\(\langle 974\rangle\)}}
\put(-1,-4.5){\makebox(0,0)[cc]{\(1;2\)}}
\put(0,-4.5){\makebox(0,0)[cc]{\(1;3\)}}
\put(-0.5,-5){\makebox(0,0)[cc]{\(\mathrm{D}.11\)}}
\put(-0.5,-5.5){\makebox(0,0)[cc]{\((243,3)\)}}

\put(1.7,-4){\makebox(0,0)[lc]{\(1;1\)}}
\put(1,-6){\makebox(0,0)[cc]{\(1;2\)}}
\put(2,-6){\makebox(0,0)[cc]{\(1;3\)}}
\put(1.5,-6.5){\makebox(0,0)[cc]{\(\mathrm{D}.11\)}}
\put(1.5,-7){\makebox(0,0)[cc]{\((729,3)\)}}

\put(3.7,-5.5){\makebox(0,0)[lc]{\(1;1\)}}
\put(3,-7.5){\makebox(0,0)[cc]{\(1;2\)}}
\put(4,-7.5){\makebox(0,0)[cc]{\(1;3\)}}
\put(3.5,-8){\makebox(0,0)[cc]{\(\mathrm{D}.11\)}}
\put(3.5,-8.5){\makebox(0,0)[cc]{\((2187,3)\)}}

\put(4.5,-4.5){\makebox(0,0)[cc]{main trunk of type \(\mathrm{b}.16\)}}
\put(5.7,-7){\makebox(0,0)[lc]{etc.}}


\put(-3,2){\makebox(0,0)[lc]{Legend:}}

\put(-1.5,2){\circle{0.2}}
\put(-1,2){\makebox(0,0)[lc]{\(\ldots\) metabelian}}

\put(-1.5,1.5){\circle{0.2}}
\put(-1.5,1.5){\circle*{0.1}}
\put(-1,1.5){\makebox(0,0)[lc]{\(\ldots\) metabelian Schur \(\sigma\)}}

\end{picture}

}

\end{figure}
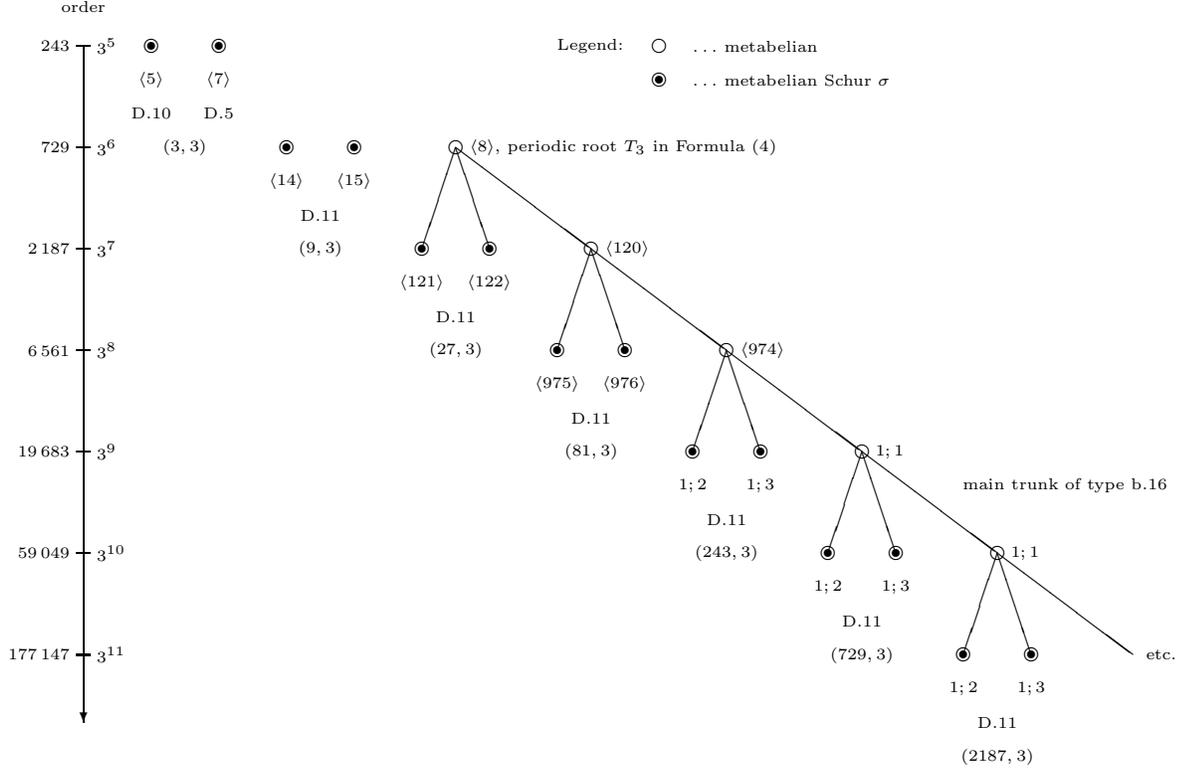

\begin{theorem}
\label{thm:b16Chain}
For each logarithmic exponent \(e\ge 3\),
the \textbf{unique} metabelian \(\mathrm{CF}\)-group \(T_e\)
with commutator quotient \(T_e/T_e^\prime\simeq (3^e,3)\hat{=}(e1)\),
punctured transfer kernel type \(\mathrm{b}.16\), \(\varkappa(T_e)\sim (004;0)\),
and logarithmic order \(\mathrm{lo}(T_e)=3+e\) 
is given by the \textbf{periodic sequence} of iterated \(p\)-descendants
\begin{equation}
\label{eqn:b16}
T_e\simeq\mathrm{SmallGroup}(729,8)(-\#1;1)^{e-3}.
\end{equation}
These groups have
first abelian quotient invariants \(\alpha_1(T_e)\sim (e1,e1,e11;(e-1)11)\).
They form the infinite main trunk of a descendant tree
with doublets as finite twigs.
For each
integer
\(e\ge 3\),
the doublet
\begin{equation}
\label{eqn:D11}
G_{e,i}\simeq\mathrm{SmallGroup}(729,8)(-\#1;1)^{e-3}-\#1;i, \qquad i\in\lbrace 2,3\rbrace,
\end{equation}
is the \textbf{unique pair} of \textbf{metabelian} Schur \(\sigma\)-groups \(G_{e,i}\) with
commutator quotient \(G_{e,i}/G_{e,i}^\prime\simeq (3^e,3)\hat{=}(e1)\).
It has punctured transfer kernel type \(\mathrm{D}.11\), \(\varkappa(G_{e,i})\sim (124;1)\),
\(\mathrm{lo}(G_{e,i})=4+e\), and
first abelian quotient invariants \(\alpha_1(G_{e,i})\sim ((e+1)1,(e+1)1,e11;e11)\).
See Figure
\ref{fig:SchurSigma3}.
\end{theorem}


\begin{proof}
Let \(s_2=\lbrack y,x\rbrack\) denote the main commutator,
\(s_3=\lbrack s_2,x\rbrack\) and \(t_3=\lbrack s_2,y\rbrack\) higher commutators,
and \(w=x^{3^e}\) the last non-trivial power.
For each \(e\ge 3\), we have parametrized presentations
\begin{equation}
\label{eqn:Pres3}
\begin{aligned}
T_{e+1} &= \langle x,y,s_2,s_3,t_3,w\mid x^{3^{e+1}}=w^3=1,\ y^3=s_3,\ t_3=s_3\rangle, \\
G_{e,i} &= \langle x,y,s_2,s_3,t_3,w\mid x^{3^{e+1}}=w^3=1,\ y^3=s_3,\ t_3=s_3w^{i-1}\rangle, \qquad i\in\lbrace 2,3\rbrace.
\end{aligned}
\end{equation}
The \(p\)-class is given by
\(\mathrm{cl}_p(G_{e,i})=e+1\) respectively \(\mathrm{cl}_p(T_{e+1})=e+1\).
For \(e\ge 3\), the last non-trivial lower \(p\)-central is given by
\(P_{e}(G_{e,i})=\langle w\rangle\) respectively \(P_{e}(T_{e+1})=\langle w\rangle\),
whence all three groups share the common \(p\)-parent
\(G_{e,i}/P_{e}(G_{e,i})\simeq T_{e+1}/P_{e}(T_{e+1})\simeq T_e\).

Repeated recursive applications of the \(p\)-group generation algorithm
eventually produce Figure
\ref{fig:SchurSigma3},
and thus confirm Formula
\eqref{eqn:b16}
and
\eqref{eqn:D11}.
All groups have \(\varrho\sim (2,2,3;3)\).
Cf.
\cite[Thm. 4--7]{Ma2021a}.
\end{proof}


\medskip
In particular, for \(e=3\) and \(e=4\), we have the pairs \\
\(\mathrm{SmallGroup}(2187,121)\simeq G_{3,2}=\langle x,y,s_2,s_3,t_3,w\mid x^{27}=w,y^3=s_3,t_3=s_3\cdot w\rangle\), \\
\(\mathrm{SmallGroup}(2187,122)\simeq G_{3,3}=\langle x,y,s_2,s_3,t_3,w\mid x^{27}=w,y^3=s_3,t_3=s_3\cdot w^2\rangle\); and \\
\(\mathrm{SmallGroup}(6561,975)\simeq G_{4,2}=\langle x,y,s_2,s_3,t_3,w\mid x^{81}=w,y^3=s_3,t_3=s_3\cdot w\rangle\), \\
\(\mathrm{SmallGroup}(6561,976)\simeq G_{4,3}=\langle x,y,s_2,s_3,t_3,w\mid x^{81}=w,y^3=s_3,t_3=s_3\cdot w^2\rangle\).

\medskip
The smallest pair, for \(e=2\), however, is exceptional: \\
\(\mathrm{SmallGroup}(729,14)\simeq\langle x,y,s_2,s_3,t_3\mid x^9=s_3,y^3=t_3\rangle\), \\
\(\mathrm{SmallGroup}(729,15)\simeq\langle x,y,s_2,s_3,t_3\mid x^9=s_3,y^3=t_3^2\rangle\).

\medskip
Theorem
\ref{thm:b16Chain}
is illustrated for \(3\le e\le 7\) by Figure
\ref{fig:SchurSigma3},
using the notation of
\cite{BEO2005,GNO2006}.


\section{Moderate Schur \(\sigma\)-groups \(G\) with \(G/G^\prime\simeq (3^e,3)\), \(e\ge 5\), ground state}
\label{s:NonMetabelian5}

\begin{figure}[hb]
\caption{Schur \(\sigma\)-groups \(G\) with \(\varrho(G)\sim (2,2,3;3)\), \(G/G^\prime\simeq (3^e,3)\), \(2\le e\le 7\)}
\label{fig:SchurSigma5}

{\tiny

\setlength{\unitlength}{0.9cm}
\begin{picture}(16,20)(-11,-17)

\put(-10,2.5){\makebox(0,0)[cb]{order}}

\put(-10,2){\line(0,-1){17}}
\multiput(-10.1,2)(0,-1.5){13}{\line(1,0){0.2}}

\put(-10.2,2){\makebox(0,0)[rc]{\(9\)}}
\put(-9.8,2){\makebox(0,0)[lc]{\(3^2\)}}
\put(-10.2,0.5){\makebox(0,0)[rc]{\(27\)}}
\put(-9.8,0.5){\makebox(0,0)[lc]{\(3^3\)}}
\put(-10.2,-1){\makebox(0,0)[rc]{\(81\)}}
\put(-9.8,-1){\makebox(0,0)[lc]{\(3^4\)}}
\put(-10.2,-2.5){\makebox(0,0)[rc]{\(243\)}}
\put(-9.8,-2.5){\makebox(0,0)[lc]{\(3^5\)}}
\put(-10.2,-4){\makebox(0,0)[rc]{\(729\)}}
\put(-9.8,-4){\makebox(0,0)[lc]{\(3^6\)}}
\put(-10.2,-5.5){\makebox(0,0)[rc]{\(2\,187\)}}
\put(-9.8,-5.5){\makebox(0,0)[lc]{\(3^7\)}}
\put(-10.2,-7){\makebox(0,0)[rc]{\(6\,561\)}}
\put(-9.8,-7){\makebox(0,0)[lc]{\(3^8\)}}
\put(-10.2,-8.5){\makebox(0,0)[rc]{\(19\,683\)}}
\put(-9.8,-8.5){\makebox(0,0)[lc]{\(3^9\)}}
\put(-10.2,-10){\makebox(0,0)[rc]{\(59\,049\)}}
\put(-9.8,-10){\makebox(0,0)[lc]{\(3^{10}\)}}
\put(-10.2,-11.5){\makebox(0,0)[rc]{\(177\,147\)}}
\put(-9.8,-11.5){\makebox(0,0)[lc]{\(3^{11}\)}}
\put(-10.2,-13){\makebox(0,0)[rc]{\(531\,441\)}}
\put(-9.8,-13){\makebox(0,0)[lc]{\(3^{12}\)}}
\put(-10.2,-14.5){\makebox(0,0)[rc]{\(1\,594\,323\)}}
\put(-9.8,-14.5){\makebox(0,0)[lc]{\(3^{13}\)}}
\put(-10.2,-16){\makebox(0,0)[rc]{\(4\,782\,969\)}}
\put(-9.8,-16){\makebox(0,0)[lc]{\(3^{14}\)}}

\put(-10,-15){\vector(0,-1){2}}

\put(-9.1,1.9){\framebox(0.2,0.2){}}
\put(-9,2){\circle*{0.1}}

\put(-9,0.5){\circle{0.2}}
\put(-9,-2.5){\circle{0.2}}
\put(-9,-4){\circle*{0.2}}
\put(-8.5,-5.5){\circle{0.2}}

\put(-7,-1){\circle{0.2}}
\put(-7,-4){\circle{0.2}}
\put(-7,-5.5){\circle*{0.2}}
\put(-6.5,-7){\circle{0.2}}

\put(-5,-4){\circle{0.2}}
\put(-5,-7){\circle*{0.2}}
\put(-4.5,-8.5){\circle{0.2}}

\put(-3,-7){\circle{0.2}}
\put(-3,-10){\circle{0.2}}

\put(-1,-10){\circle{0.2}}
\put(-1,-11.5){\circle{0.2}}

\put(1,-11.5){\circle{0.2}}
\put(1,-13){\circle{0.2}}

\put(3,-13){\circle{0.2}}
\put(3,-14.5){\circle{0.2}}

\put(-9.1,-7.1){\framebox(0.2,0.2){}}
\put(-7.1,-8.6){\framebox(0.2,0.2){}}
\put(-5.1,-10.1){\framebox(0.2,0.2){}}
\put(-3.1,-11.6){\framebox(0.2,0.2){}}
\put(-1.1,-13.1){\framebox(0.2,0.2){}}
\put(0.9,-14.6){\framebox(0.2,0.2){}}
\put(2.9,-16.1){\framebox(0.2,0.2){}}


\put(-9,2){\line(0,-1){1.5}}
\put(-9,0.5){\line(0,-1){3}}
\put(-9,-2.5){\line(0,-1){1.5}}
\put(-9,-4){\line(0,-1){3}}
\put(-9,-4){\line(1,-3){0.5}}

\put(-9,2){\line(2,-3){2}}
\put(-7,-1){\line(0,-1){3}}
\put(-7,-4){\line(0,-1){1.5}}
\put(-7,-5.5){\line(0,-1){3}}
\put(-7,-5.5){\line(1,-3){0.5}}

\put(-7,-1){\line(2,-3){2}}
\put(-5,-4){\line(0,-1){3}}
\put(-5,-7){\line(0,-1){3}}
\put(-5,-7){\line(1,-3){0.5}}

\put(-5,-4){\line(2,-3){2}}
\put(-3,-7){\line(0,-1){3}}
\put(-3,-10){\line(0,-1){1.5}}

\put(-3,-7){\line(2,-3){2}}
\put(-1,-10){\line(0,-1){1.5}}
\put(-1,-11.5){\line(0,-1){1.5}}

\put(-1,-10){\line(4,-3){2}}
\put(1,-11.5){\line(0,-1){1.5}}
\put(1,-13){\line(0,-1){1.5}}

\put(1,-11.5){\line(4,-3){2}}
\put(3,-13){\line(0,-1){1.5}}
\put(3,-14.5){\line(0,-1){1.5}}

\put(3,-13){\line(4,-3){1}}


\put(-8.7,2){\makebox(0,0)[lc]{\(\langle 2\rangle\)}}

\put(-8.7,0.5){\makebox(0,0)[lc]{\(\langle 3\rangle\)}}
\put(-8.7,-2.5){\makebox(0,0)[lc]{\(\langle 6\rangle\)}}
\put(-8.7,-4){\makebox(0,0)[lc]{\(\langle 49\rangle\)}}
\put(-8.2,-5.5){\makebox(0,0)[lc]{\(\langle 289\vert\)}}
\put(-8.1,-5.8){\makebox(0,0)[lc]{\(290\rangle\)}}

\put(-6.7,-1){\makebox(0,0)[lc]{\(\langle 3\rangle\)}}
\put(-6.7,-4){\makebox(0,0)[lc]{\(\langle 13\rangle\)}}
\put(-6.7,-5.5){\makebox(0,0)[lc]{\(\langle 168\rangle\)}}
\put(-6.2,-7){\makebox(0,0)[lc]{\(\langle 1689\vert\)}}
\put(-6.1,-7.3){\makebox(0,0)[lc]{\(1690\rangle\)}}

\put(-4.7,-4){\makebox(0,0)[lc]{\(\langle 7\rangle\)}}
\put(-4.7,-7){\makebox(0,0)[lc]{\(\langle 98\rangle\)}}
\put(-4.2,-8.5){\makebox(0,0)[lc]{\(1;4\vert 5\)}}

\put(-2.7,-7){\makebox(0,0)[lc]{\(\langle 93\rangle\)}}
\put(-2.7,-10){\makebox(0,0)[lc]{\(2;7\vert 8\)}}

\put(-0.7,-10){\makebox(0,0)[lc]{\(2;2\), periodic root \(T_5\) in Formula \eqref{eqn:b16GroundD10}}}
\put(-0.7,-11.5){\makebox(0,0)[lc]{\(1;2\vert 3\)}}

\put(1.3,-11.5){\makebox(0,0)[lc]{\(1;1\)}}
\put(1.3,-13){\makebox(0,0)[lc]{\(1;2\vert 3\)}}

\put(3.3,-13){\makebox(0,0)[lc]{\(1;1\)}}
\put(3.3,-14.5){\makebox(0,0)[lc]{\(1;2\vert 3\)}}

\put(-8.7,-7){\makebox(0,0)[lc]{\(\langle 617\vert 618\rangle\)}}
\put(-9,-7.5){\makebox(0,0)[lc]{\(\mathrm{E}.14\)}}
\put(-9,-8){\makebox(0,0)[lc]{\((3,3)\)}}

\put(-6.7,-8.5){\makebox(0,0)[lc]{\(2;8\vert 9\)}}
\put(-7,-9){\makebox(0,0)[lc]{\(\mathrm{D}.10\)}}
\put(-7,-9.5){\makebox(0,0)[lc]{\((9,3)\)}}

\put(-4.7,-10){\makebox(0,0)[lc]{\(2;2\vert 3\)}}
\put(-5,-10.5){\makebox(0,0)[lc]{\(\mathrm{D}.10\)}}
\put(-5,-11){\makebox(0,0)[lc]{\((27,3)\)}}

\put(-2.7,-11.5){\makebox(0,0)[lc]{\(1;1\)}}
\put(-3,-12){\makebox(0,0)[lc]{\(\mathrm{D}.10\)}}
\put(-3,-12.5){\makebox(0,0)[lc]{\((81,3)\)}}

\put(-0.7,-13){\makebox(0,0)[lc]{\(1;1\)}}
\put(-1,-13.5){\makebox(0,0)[lc]{\(\mathrm{D}.10\)}}
\put(-1,-14){\makebox(0,0)[lc]{\((243,3)\)}}

\put(1.3,-14.5){\makebox(0,0)[lc]{\(1;1\)}}
\put(1,-15){\makebox(0,0)[lc]{\(\mathrm{D}.10\)}}
\put(1,-15.5){\makebox(0,0)[lc]{\((729,3)\)}}

\put(3.3,-16){\makebox(0,0)[lc]{\(1;1\)}}
\put(3,-16.5){\makebox(0,0)[lc]{\(\mathrm{D}.10\)}}
\put(3,-17){\makebox(0,0)[lc]{\((2187,3)\)}}

\put(3.5,-12){\makebox(0,0)[cc]{main trunk of type \(\mathrm{b}.16\)}}
\put(4.3,-13.7){\makebox(0,0)[lt]{etc.}}

\put(-3,-0.6){\makebox(0,0)[lc]{Legend:}}

\put(-1.6,-0.7){\framebox(0.2,0.2){}}
\put(-1.5,-0.6){\circle*{0.1}}
\put(-1,-0.6){\makebox(0,0)[lc]{\(\ldots\) abelian}}

\put(-1.5,-1){\circle{0.2}}
\put(-1,-1){\makebox(0,0)[lc]{\(\ldots\) metabelian}}

\put(-1.5,-1.4){\circle*{0.2}}
\put(-1,-1.4){\makebox(0,0)[lc]{\(\ldots\) metabelian with bifurcation}}

\put(-1.6,-1.9){\framebox(0.2,0.2){}}
\put(-1,-1.8){\makebox(0,0)[lc]{\(\ldots\) non-metabelian Schur \(\sigma\)}}

\end{picture}

}

\end{figure}
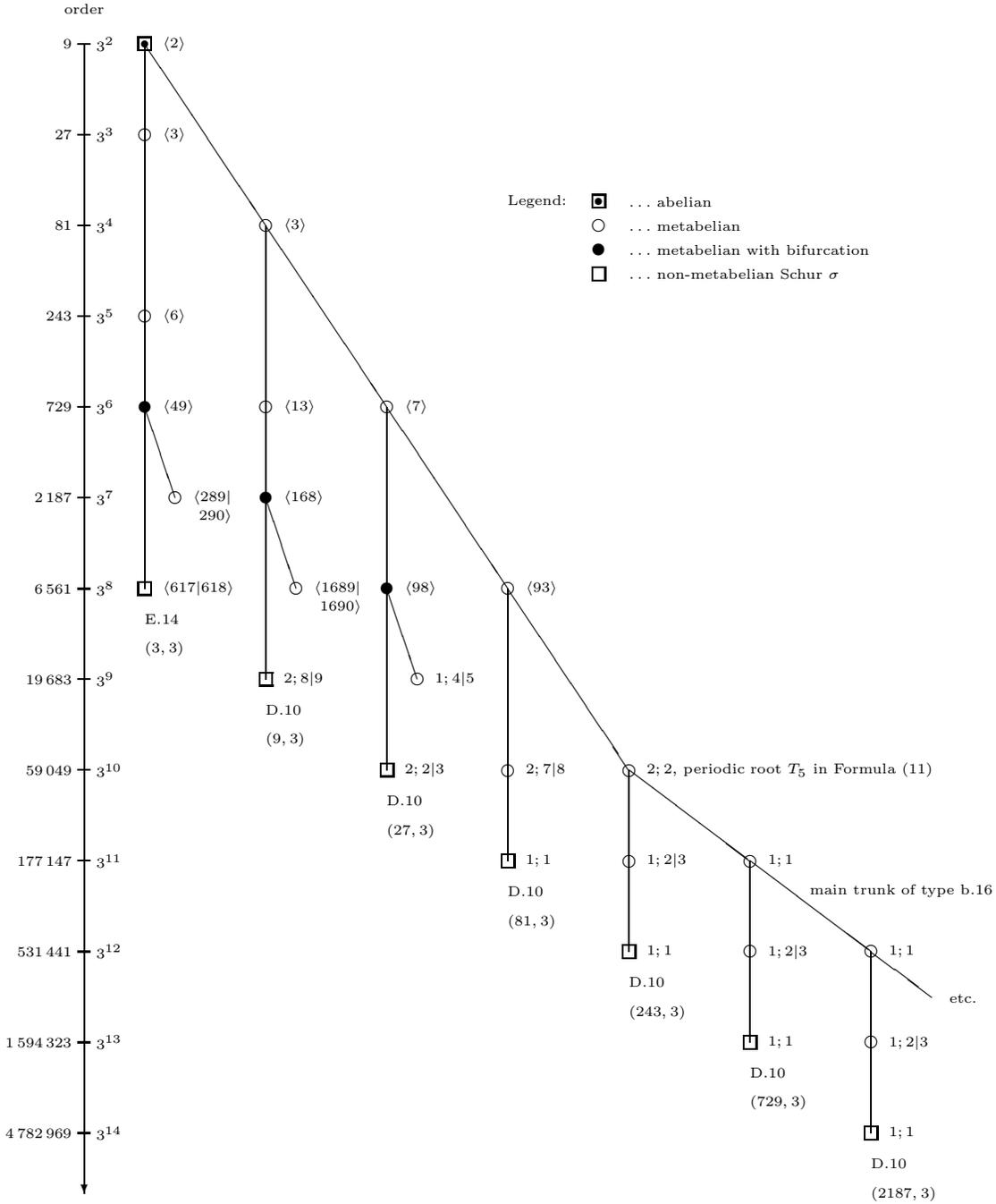

\begin{figure}[hb]
\caption{Schur \(\sigma\)-groups \(G\) with \(\varrho(G)\sim (2,2,2;3)\), \(G/G^\prime\simeq (3^e,3)\), \(2\le e\le 7\)}
\label{fig:SchurSigma5D6}

{\tiny

\setlength{\unitlength}{0.9cm}
\begin{picture}(16,20)(-11,-17)

\put(-10,2.5){\makebox(0,0)[cb]{order}}

\put(-10,2){\line(0,-1){17}}
\multiput(-10.1,2)(0,-1.5){13}{\line(1,0){0.2}}

\put(-10.2,2){\makebox(0,0)[rc]{\(9\)}}
\put(-9.8,2){\makebox(0,0)[lc]{\(3^2\)}}
\put(-10.2,0.5){\makebox(0,0)[rc]{\(27\)}}
\put(-9.8,0.5){\makebox(0,0)[lc]{\(3^3\)}}
\put(-10.2,-1){\makebox(0,0)[rc]{\(81\)}}
\put(-9.8,-1){\makebox(0,0)[lc]{\(3^4\)}}
\put(-10.2,-2.5){\makebox(0,0)[rc]{\(243\)}}
\put(-9.8,-2.5){\makebox(0,0)[lc]{\(3^5\)}}
\put(-10.2,-4){\makebox(0,0)[rc]{\(729\)}}
\put(-9.8,-4){\makebox(0,0)[lc]{\(3^6\)}}
\put(-10.2,-5.5){\makebox(0,0)[rc]{\(2\,187\)}}
\put(-9.8,-5.5){\makebox(0,0)[lc]{\(3^7\)}}
\put(-10.2,-7){\makebox(0,0)[rc]{\(6\,561\)}}
\put(-9.8,-7){\makebox(0,0)[lc]{\(3^8\)}}
\put(-10.2,-8.5){\makebox(0,0)[rc]{\(19\,683\)}}
\put(-9.8,-8.5){\makebox(0,0)[lc]{\(3^9\)}}
\put(-10.2,-10){\makebox(0,0)[rc]{\(59\,049\)}}
\put(-9.8,-10){\makebox(0,0)[lc]{\(3^{10}\)}}
\put(-10.2,-11.5){\makebox(0,0)[rc]{\(177\,147\)}}
\put(-9.8,-11.5){\makebox(0,0)[lc]{\(3^{11}\)}}
\put(-10.2,-13){\makebox(0,0)[rc]{\(531\,441\)}}
\put(-9.8,-13){\makebox(0,0)[lc]{\(3^{12}\)}}
\put(-10.2,-14.5){\makebox(0,0)[rc]{\(1\,594\,323\)}}
\put(-9.8,-14.5){\makebox(0,0)[lc]{\(3^{13}\)}}
\put(-10.2,-16){\makebox(0,0)[rc]{\(4\,782\,969\)}}
\put(-9.8,-16){\makebox(0,0)[lc]{\(3^{14}\)}}

\put(-10,-15){\vector(0,-1){2}}

\put(-9.1,1.9){\framebox(0.2,0.2){}}
\put(-9,2){\circle*{0.1}}

\put(-9,0.5){\circle{0.2}}
\put(-9,-2.5){\circle{0.2}}
\put(-9,-4){\circle*{0.2}}
\put(-8.5,-5.5){\circle{0.2}}

\put(-7,-1){\circle{0.2}}
\put(-7,-4){\circle{0.2}}
\put(-7,-5.5){\circle*{0.2}}
\put(-6.5,-7){\circle{0.2}}

\put(-5,-4){\circle{0.2}}
\put(-5,-7){\circle*{0.2}}
\put(-4.5,-8.5){\circle{0.2}}

\put(-3,-7){\circle{0.2}}
\put(-3,-10){\circle{0.2}}

\put(-1,-10){\circle{0.2}}
\put(-1,-11.5){\circle{0.2}}

\put(1,-11.5){\circle{0.2}}
\put(1,-13){\circle{0.2}}

\put(3,-13){\circle{0.2}}
\put(3,-14.5){\circle{0.2}}

\put(-9.1,-7.1){\framebox(0.2,0.2){}}
\put(-7.1,-8.6){\framebox(0.2,0.2){}}
\put(-5.1,-10.1){\framebox(0.2,0.2){}}
\put(-3.1,-11.6){\framebox(0.2,0.2){}}
\put(-1.1,-13.1){\framebox(0.2,0.2){}}
\put(0.9,-14.6){\framebox(0.2,0.2){}}
\put(2.9,-16.1){\framebox(0.2,0.2){}}


\put(-9,2){\line(0,-1){1.5}}
\put(-9,0.5){\line(0,-1){3}}
\put(-9,-2.5){\line(0,-1){1.5}}
\put(-9,-4){\line(0,-1){3}}
\put(-9,-4){\line(1,-3){0.5}}

\put(-9,2){\line(2,-3){2}}
\put(-7,-1){\line(0,-1){3}}
\put(-7,-4){\line(0,-1){1.5}}
\put(-7,-5.5){\line(0,-1){3}}
\put(-7,-5.5){\line(1,-3){0.5}}

\put(-7,-1){\line(2,-3){2}}
\put(-5,-4){\line(0,-1){3}}
\put(-5,-7){\line(0,-1){3}}
\put(-5,-7){\line(1,-3){0.5}}

\put(-5,-4){\line(2,-3){2}}
\put(-3,-7){\line(0,-1){3}}
\put(-3,-10){\line(0,-1){1.5}}

\put(-3,-7){\line(2,-3){2}}
\put(-1,-10){\line(0,-1){1.5}}
\put(-1,-11.5){\line(0,-1){1.5}}

\put(-1,-10){\line(4,-3){2}}
\put(1,-11.5){\line(0,-1){1.5}}
\put(1,-13){\line(0,-1){1.5}}

\put(1,-11.5){\line(4,-3){2}}
\put(3,-13){\line(0,-1){1.5}}
\put(3,-14.5){\line(0,-1){1.5}}

\put(3,-13){\line(4,-3){1}}


\put(-8.7,2){\makebox(0,0)[lc]{\(\langle 2\rangle\)}}

\put(-8.7,0.5){\makebox(0,0)[lc]{\(\langle 3\rangle\)}}
\put(-8.7,-2.5){\makebox(0,0)[lc]{\(\langle 8\rangle\)}}
\put(-8.7,-4){\makebox(0,0)[lc]{\(\langle 54\rangle\)}}
\put(-8.2,-5.5){\makebox(0,0)[lc]{\(\langle 302\vert\)}}
\put(-8.1,-5.8){\makebox(0,0)[lc]{\(306\rangle\)}}

\put(-6.7,-1){\makebox(0,0)[lc]{\(\langle 3\rangle\)}}
\put(-6.7,-4){\makebox(0,0)[lc]{\(\langle 18\vert 21\rangle\)}}
\put(-6.7,-5.5){\makebox(0,0)[lc]{\(\langle 181\vert 191\rangle\)}}
\put(-6.2,-7){\makebox(0,0)[lc]{\(\langle 1741\vert\)}}
\put(-6.1,-7.3){\makebox(0,0)[lc]{\(1779\rangle\)}}

\put(-4.7,-4){\makebox(0,0)[lc]{\(\langle 6\rangle\)}}
\put(-4.7,-7){\makebox(0,0)[lc]{\(\langle 88\vert 91\rangle\)}}
\put(-4.2,-8.5){\makebox(0,0)[lc]{\(1;4\)}}

\put(-2.7,-7){\makebox(0,0)[lc]{\(\langle 85\rangle\)}}
\put(-2.7,-10){\makebox(0,0)[lc]{\(2;8\vert 12\)}}

\put(-0.7,-10){\makebox(0,0)[lc]{\(2;4\), periodic root \(T_5\) in Formula \eqref{eqn:a1GroundD6}}}
\put(-0.7,-11.5){\makebox(0,0)[lc]{\(1;2\vert 3\)}}

\put(1.3,-11.5){\makebox(0,0)[lc]{\(1;1\)}}
\put(1.3,-13){\makebox(0,0)[lc]{\(1;2\vert 3\)}}

\put(3.3,-13){\makebox(0,0)[lc]{\(1;1\)}}
\put(3.3,-14.5){\makebox(0,0)[lc]{\(1;2\vert 3\)}}

\put(-8.7,-7){\makebox(0,0)[lc]{\(\langle 620\vert 624\rangle\)}}
\put(-9,-7.5){\makebox(0,0)[lc]{\(\mathrm{E}.9\)}}
\put(-9,-8){\makebox(0,0)[lc]{\((3,3)\)}}

\put(-6.7,-8.5){\makebox(0,0)[lc]{\(2;4\)}}
\put(-7,-9){\makebox(0,0)[lc]{\(\mathrm{D}.6\)}}
\put(-7,-9.5){\makebox(0,0)[lc]{\((9,3)\)}}

\put(-4.7,-10){\makebox(0,0)[lc]{\(2;4\)}}
\put(-5,-10.5){\makebox(0,0)[lc]{\(\mathrm{D}.6\)}}
\put(-5,-11){\makebox(0,0)[lc]{\((27,3)\)}}

\put(-2.7,-11.5){\makebox(0,0)[lc]{\(1;1\)}}
\put(-3,-12){\makebox(0,0)[lc]{\(\mathrm{D}.6\)}}
\put(-3,-12.5){\makebox(0,0)[lc]{\((81,3)\)}}

\put(-0.7,-13){\makebox(0,0)[lc]{\(1;1\)}}
\put(-1,-13.5){\makebox(0,0)[lc]{\(\mathrm{D}.6\)}}
\put(-1,-14){\makebox(0,0)[lc]{\((243,3)\)}}

\put(1.3,-14.5){\makebox(0,0)[lc]{\(1;1\)}}
\put(1,-15){\makebox(0,0)[lc]{\(\mathrm{D}.6\)}}
\put(1,-15.5){\makebox(0,0)[lc]{\((729,3)\)}}

\put(3.3,-16){\makebox(0,0)[lc]{\(1;1\)}}
\put(3,-16.5){\makebox(0,0)[lc]{\(\mathrm{D}.6\)}}
\put(3,-17){\makebox(0,0)[lc]{\((2187,3)\)}}

\put(3.5,-12){\makebox(0,0)[cc]{main trunk of type \(\mathrm{a}.1\)}}
\put(4.3,-13.7){\makebox(0,0)[lt]{etc.}}

\put(-3,-0.6){\makebox(0,0)[lc]{Legend:}}

\put(-1.6,-0.7){\framebox(0.2,0.2){}}
\put(-1.5,-0.6){\circle*{0.1}}
\put(-1,-0.6){\makebox(0,0)[lc]{\(\ldots\) abelian}}

\put(-1.5,-1){\circle{0.2}}
\put(-1,-1){\makebox(0,0)[lc]{\(\ldots\) metabelian}}

\put(-1.5,-1.4){\circle*{0.2}}
\put(-1,-1.4){\makebox(0,0)[lc]{\(\ldots\) metabelian with bifurcation}}

\put(-1.6,-1.9){\framebox(0.2,0.2){}}
\put(-1,-1.8){\makebox(0,0)[lc]{\(\ldots\) non-metabelian Schur \(\sigma\)}}

\end{picture}

}

\end{figure}
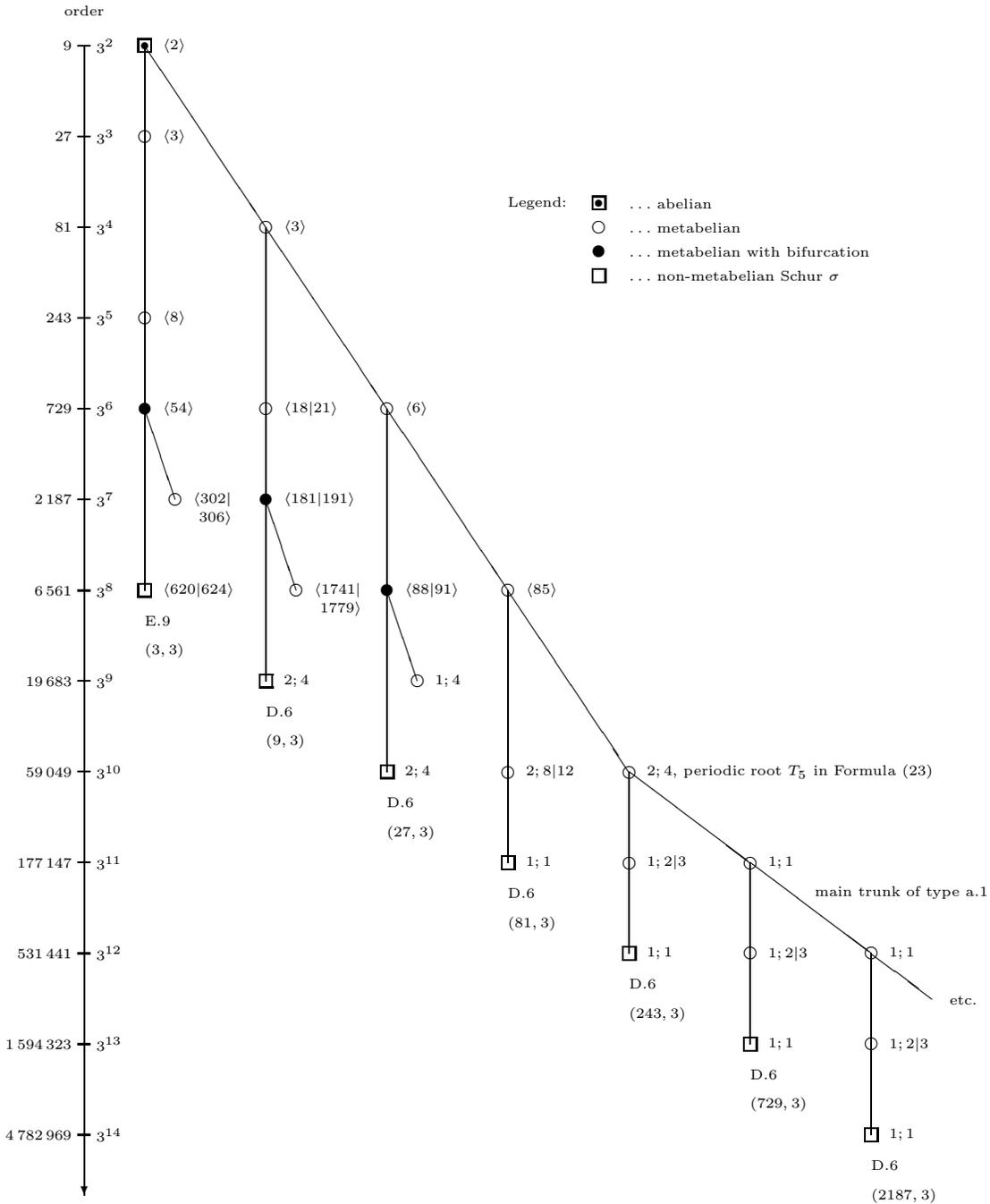

\noindent
Figure
\ref{fig:SchurSigma5},
resp.
\ref{fig:SchurSigma5D6},
shows that the construction process for the two
non-metabelian Schur \(\sigma\)-groups \(G\) with order \(\#G=3^{7+e}\) and
punctured transfer kernel type \(\mathrm{D}.10\), resp. \(\mathrm{D}.6\), 
becomes increasingly difficult for the commutator quotients \(G/G^\prime\simeq (27,3)\), \((81,3)\), \((243,3)\).
For the commutator quotient \(G/G^\prime\simeq (729,3)\), however,
an \textbf{unexpected tranquilization occurs},
and the construction process becomes settled with a \textbf{simple step size one periodicity}.


\begin{theorem}
\label{thm:Periodicity5}
The four pairs of Schur \(\sigma\)-groups \(G\) with soluble length \(\mathrm{sl}(G)=3\),
commutator quotient \(G/G^\prime\simeq (3^e,3)\hat{=}(e1)\), \(e\ge 5\), 
punctured transfer kernel types \(\mathrm{D}.10\), \(\mathrm{D}.5\), \(\mathrm{C}.4\), \(\mathrm{D}.6\),
and order \(\#G=3^{7+e}\)
are given by the following \textbf{bottom up} construction process.
\begin{itemize} 
\item
For type \(\mathrm{D}.10\), \(\varkappa\sim (411;3)\),
\(\alpha\sim (e22,(e+1)1,(e+1)1;e11)\), \(\varrho\sim (3,2,2;3)\), \\
let \(A_2:=\mathrm{SmallGroup}(6561,93)-\#2;2\), then
\begin{equation}
\label{eqn:D10Clone}
G\simeq A_2(-\#1;1)^{e-5}-\#1;i-\#1;1, \qquad i\in\lbrace 2,3\rbrace.
\end{equation}
\item
For type \(\mathrm{D}.5\), \(\varkappa\sim (211;3)\), 
\(\alpha\sim ((e+1)21,(e+1)1,(e+1)1;e11)\), \(\varrho\sim (3,2,2;3)\), \\
let \(A_4:=\mathrm{SmallGroup}(6561,93)-\#2;4\), then
\begin{equation}
\label{eqn:D5Clone}
G\simeq A_4(-\#1;1)^{e-5}-\#1;i-\#1;1, \qquad i\in\lbrace 2,3\rbrace.
\end{equation}
\item
For type \(\mathrm{C}.4\), \(\varkappa\sim(311;3)\),
\(\alpha\sim ((e+1)21,(e+1)1,(e+1)1;e11)\), \(\varrho\sim (3,2,2;3)\), \\
let \(A_5:=\mathrm{SmallGroup}(6561,93)-\#2;5\), then
\begin{equation}
\label{eqn:C4Clone}
G\simeq A_5(-\#1;1)^{e-5}-\#1;i-\#1;1, \qquad i\in\lbrace 2,3\rbrace.
\end{equation}
\item
For type \(\mathrm{D}.6\), \(\varkappa\sim (123;1)\), 
\(\alpha\sim ((e+1)1,(e+1)1,(e+1)1;e22)\), \(\varrho\sim (2,2,2;3)\), \\
let \(A_0:=\mathrm{SmallGroup}(6561,85)-\#2;4\), then
\begin{equation}
\label{eqn:D6Clone}
G\simeq A_0(-\#1;1)^{e-5}-\#1;i-\#1;1, \qquad i\in\lbrace 2,3\rbrace.
\end{equation}
\end{itemize} 
\end{theorem}


\begin{proof}
(Theorem
\ref{thm:Periodicity5})
The proof consists of the construction of successive descendants of
\(A_2\), \(A_4\), \(A_5\), \(A_0\)
in the way indicated in Theorem
\ref{thm:Periodicity5}
by means of the \(p\)-group generation algorithm
\cite{HEO2005}
by Newman
\cite{Nm1977}
and O'Brien
\cite{Ob1990},
implemented in Magma
\cite{BCP1997,BCFS2021,MAGMA2021}.
Cf.
\cite[Thm. 12--13]{Ma2021a}.
\end{proof}


After Theorem
\ref{thm:Periodicity5},
we begin with a thorough investigation of the ground state of type \(\mathrm{D}.10\).

\begin{theorem}
\label{thm:b16ChainGroundD10}
For each logarithmic exponent \(e\ge 5\),
a metabelian \(\mathrm{CF}\)-group \(T_e\)
with commutator quotient \(T_e/T_e^\prime\simeq (3^e,3)\hat{=}(e1)\),
punctured transfer kernel type \(\mathrm{b}.16\), \(\varkappa(T_e)\sim (004;0)\),
and rank distribution \(\varrho\sim (2,2,3;3)\)
is given by the \textbf{periodic sequence} of iterated \(p\)-descendants
\begin{equation}
\label{eqn:b16GroundD10}
T_e\simeq\mathrm{SmallGroup}(6561,93)-\#2;2(-\#1;1)^{e-5}.
\end{equation}
These groups have logarithmic order \(\mathrm{lo}(T_e)=5+e\) and
first abelian quotient invariants \(\alpha_1(T_e)\sim (e1,e1,e22;(e-1)11)\).
They form the infinite main trunk of a descendant tree
with finite double-twigs of depth two.
For each integer \(e\ge 5\), the doublet
\begin{equation}
\label{eqn:GroundD10M}
M_{e,i}\simeq\mathrm{SmallGroup}(6561,93)-\#2;2(-\#1;1)^{e-5}-\#1;i, \qquad i\in\lbrace 2,3\rbrace,
\end{equation}
respectively
\begin{equation}
\label{eqn:GroundD10}
G_{e,i}\simeq\mathrm{SmallGroup}(6561,93)-\#2;2(-\#1;1)^{e-5}-\#1;i-\#1;1, \qquad i\in\lbrace 2,3\rbrace,
\end{equation}
is the \textbf{unique pair} of \textbf{metabelian} \(\sigma\)-groups \(M_{e,i}\simeq G_{e,i}/G_{e,i}^{\prime\prime}\),
respectively \textbf{non-metabelian} Schur \(\sigma\)-groups \(G_{e,i}\), with
commutator quotient \((3^e,3)\hat{=}(e1)\),
punctured transfer kernel type \(\mathrm{D}.10\), \(\varkappa\sim (114;3)\),
\(\mathrm{lo}(M_{e,i})=6+e\), respectively \(\mathrm{lo}(G_{e,i})=7+e\), and
first abelian quotient invariants \(\alpha_1\sim ((e+1)1,(e+1)1,e22;e11)\).
See Figure
\ref{fig:SchurSigma5}.
\end{theorem}


\begin{proof}
Let \(s_2=\lbrack y,x\rbrack\) denote the main commutator,
\(s_3=\lbrack s_2,x\rbrack\), \(s_4=\lbrack s_3,x\rbrack\), \(s_5=\lbrack s_4,x\rbrack\)
and \(t_3=\lbrack s_2,y\rbrack\), \(t_4=\lbrack s_3,y\rbrack\), \(t_5=\lbrack s_4,y\rbrack\) higher commutators,
and \(w=x^{3^e}\) the last non-trivial power.
For each \(e\ge 5\), we have parametrized presentations with \(i\in\lbrace 2,3\rbrace\)
\begin{equation}
\label{eqn:Pres5D10}
\begin{aligned}
T_{e+1} &= \langle x,y\mid w^3=1,\ y^3=s_5,\ s_2^3=s_4^2s_5,\ s_3^3=s_5^2,\ t_3=s_3,\ t_4=s_4,\ t_5=s_5\rangle, \\
M_{e,i} &= \langle x,y\mid w^3=1,\ y^3=s_5,\ s_2^3=s_4^2s_5,\ s_3^3=s_5^2,\ t_3=s_3w^{i-1},\ t_4=s_4,\ t_5=s_5\rangle.
\end{aligned}
\end{equation}
The \(p\)-class is given by
\(\mathrm{cl}_p(M_{e,i})=e+1\) respectively \(\mathrm{cl}_p(T_{e+1})=e+1\).
For \(e\ge 5\), the last non-trivial lower \(p\)-central is given by
\(P_{e}(M_{e,i})=\langle w\rangle\) respectively \(P_{e}(T_{e+1})=\langle w\rangle\),
whence all three groups share the common \(p\)-parent
\(M_{e,i}/P_{e}(M_{e,i})\simeq T_{e+1}/P_{e}(T_{e+1})\simeq T_e\).

Repeated recursive applications of the \(p\)-group generation algorithm
eventually produce Figure
\ref{fig:SchurSigma5},
and thus confirm Formulas
\eqref{eqn:b16GroundD10},
\eqref{eqn:GroundD10M}
and
\eqref{eqn:GroundD10}.
All groups have \(\varrho\sim (2,2,3;3)\).
Cf.
\cite[Thm. 11--13]{Ma2021a}.
\end{proof}


We continue with the ground state of type \(\mathrm{C}.4\).

\begin{theorem}
\label{thm:a1ChainGroundC4}
For each logarithmic exponent \(e\ge 5\),
a metabelian \(\mathrm{CF}\)-group \(T_e\)
with commutator quotient \(T_e/T_e^\prime\simeq (3^e,3)\hat{=}(e1)\),
punctured transfer kernel type \(\mathrm{a}.1\), \(\varkappa(T_e)\sim (000;0)\),
and rank distribution \(\varrho\sim (2,2,3;3)\)
is given by the \textbf{periodic sequence} of iterated \(p\)-descendants
\begin{equation}
\label{eqn:a1GroundC4}
T_e\simeq\mathrm{SmallGroup}(6561,93)-\#2;5(-\#1;1)^{e-5}.
\end{equation}
These groups have logarithmic order \(\mathrm{lo}(T_e)=5+e\) and
first abelian quotient invariants \(\alpha_1(T_e)\sim (e1,e1,e21;(e-1)11)\).
They form the infinite main trunk of a descendant tree
with finite double-twigs of depth two.
For each integer \(e\ge 5\), the doublet
\begin{equation}
\label{eqn:GroundC4M}
M_{e,i}\simeq\mathrm{SmallGroup}(6561,93)-\#2;5(-\#1;1)^{e-5}-\#1;i, \qquad i\in\lbrace 2,3\rbrace,
\end{equation}
respectively
\begin{equation}
\label{eqn:GroundC4}
G_{e,i}\simeq\mathrm{SmallGroup}(6561,93)-\#2;5(-\#1;1)^{e-5}-\#1;i-\#1;1, \qquad i\in\lbrace 2,3\rbrace,
\end{equation}
is the \textbf{unique pair} of \textbf{metabelian} \(\sigma\)-groups \(M_{e,i}\simeq G_{e,i}/G_{e,i}^{\prime\prime}\),
respectively \textbf{non-metabelian} Schur \(\sigma\)-groups \(G_{e,i}\), with
commutator quotient \((3^e,3)\hat{=}(e1)\),
punctured transfer kernel type \(\mathrm{C}.4\), \(\varkappa\sim (113;3)\),
\(\mathrm{lo}(M_{e,i})=6+e\), respectively \(\mathrm{lo}(G_{e,i})=7+e\), and
first abelian quotient invariants \(\alpha_1\sim ((e+1)1,(e+1)1,(e+1)21;e11)\).
\end{theorem}


\begin{proof}
Let \(s_2=\lbrack y,x\rbrack\) denote the main commutator,
\(s_3=\lbrack s_2,x\rbrack\), \(s_4=\lbrack s_3,x\rbrack\), \(s_5=\lbrack s_4,x\rbrack\)
and \(t_3=\lbrack s_2,y\rbrack\), \(t_4=\lbrack s_3,y\rbrack\), \(t_5=\lbrack s_4,y\rbrack\) higher commutators,
and \(w=x^{3^e}\) the last non-trivial power.
For each \(e\ge 5\), we have parametrized presentations with \(i\in\lbrace 2,3\rbrace\)
\begin{equation}
\label{eqn:Pres5C4}
\begin{aligned}
T_{e+1} &= \langle x,y\mid w^3=1,\ y^3=s_5,\ s_2^3=s_4^2s_5,\ s_3^3=s_5^2,\ t_3=s_3s_5,\ t_4=s_4,\ t_5=s_5\rangle, \\
M_{e,i} &= \langle x,y\mid w^3=1,\ y^3=s_5,\ s_2^3=s_4^2s_5,\ s_3^3=s_5^2,\ t_3=s_3s_5w^{i-1},\ t_4=s_4,\ t_5=s_5\rangle.
\end{aligned}
\end{equation}
The \(p\)-class is given by
\(\mathrm{cl}_p(M_{e,i})=e+1\) respectively \(\mathrm{cl}_p(T_{e+1})=e+1\).
For \(e\ge 5\), the last non-trivial lower \(p\)-central is given by
\(P_{e}(M_{e,i})=\langle w\rangle\) respectively \(P_{e}(T_{e+1})=\langle w\rangle\),
whence all three groups share the common \(p\)-parent
\(M_{e,i}/P_{e}(M_{e,i})\simeq T_{e+1}/P_{e}(T_{e+1})\simeq T_e\).

Repeated recursive applications of the \(p\)-group generation algorithm
eventually produce a graph isomorphic to the tree in Figure
\ref{fig:SchurSigma5},
and thus confirm Formulas
\eqref{eqn:a1GroundC4},
\eqref{eqn:GroundC4M}
and
\eqref{eqn:GroundC4}.
All groups have \(\varrho\sim (2,2,3;3)\).
Cf.
\cite[Thm. 11--13]{Ma2021a}.
\end{proof}


Next, we look at the ground state of type \(\mathrm{D}.5\).

\begin{theorem}
\label{thm:a1ChainGroundD5}
For each logarithmic exponent \(e\ge 5\),
a metabelian \(\mathrm{CF}\)-group \(T_e\)
with commutator quotient \(T_e/T_e^\prime\simeq (3^e,3)\hat{=}(e1)\),
punctured transfer kernel type \(\mathrm{a}.1\), \(\varkappa(T_e)\sim (000;0)\),
and rank distribution \(\varrho\sim (2,2,3;3)\)
is given by the \textbf{periodic sequence} of iterated \(p\)-descendants
\begin{equation}
\label{eqn:a1GroundD5}
T_e\simeq\mathrm{SmallGroup}(6561,93)-\#2;4(-\#1;1)^{e-5}.
\end{equation}
These groups have logarithmic order \(\mathrm{lo}(T_e)=5+e\) and
first abelian quotient invariants \(\alpha_1(T_e)\sim (e1,e1,e21;(e-1)11)\).
They form the infinite main trunk of a descendant tree
with finite double-twigs of depth two.
For each integer \(e\ge 5\), the doublet
\begin{equation}
\label{eqn:GroundD5M}
M_{e,i}\simeq\mathrm{SmallGroup}(6561,93)-\#2;4(-\#1;1)^{e-5}-\#1;i, \qquad i\in\lbrace 2,3\rbrace,
\end{equation}
respectively
\begin{equation}
\label{eqn:GroundD5}
G_{e,i}\simeq\mathrm{SmallGroup}(6561,93)-\#2;4(-\#1;1)^{e-5}-\#1;i-\#1;1, \qquad i\in\lbrace 2,3\rbrace,
\end{equation}
is the \textbf{unique pair} of \textbf{metabelian} \(\sigma\)-groups \(M_{e,i}\simeq G_{e,i}/G_{e,i}^{\prime\prime}\),
respectively \textbf{non-metabelian} Schur \(\sigma\)-groups \(G_{e,i}\), with
commutator quotient \((3^e,3)\hat{=}(e1)\),
punctured transfer kernel type \(\mathrm{D}.5\), \(\varkappa\sim (112;3)\),
\(\mathrm{lo}(M_{e,i})=6+e\), respectively \(\mathrm{lo}(G_{e,i})=7+e\), and
first abelian quotient invariants \(\alpha_1\sim ((e+1)1,(e+1)1,(e+1)21;e11)\).
\end{theorem}


\begin{proof}
Let \(s_2=\lbrack y,x\rbrack\) denote the main commutator,
\(s_3=\lbrack s_2,x\rbrack\), \(s_4=\lbrack s_3,x\rbrack\), \(s_5=\lbrack s_4,x\rbrack\)
and \(t_3=\lbrack s_2,y\rbrack\), \(t_4=\lbrack s_3,y\rbrack\), \(t_5=\lbrack s_4,y\rbrack\) higher commutators,
and \(w=x^{3^e}\) the last non-trivial power.
For each \(e\ge 5\), we have parametrized presentations with \(i\in\lbrace 2,3\rbrace\)
\begin{equation}
\label{eqn:Pres5D5}
\begin{aligned}
T_{e+1} &= \langle x,y\mid w^3=1,\ y^3=s_5^2,\ s_2^3=s_4^2s_5,\ s_3^3=s_5^2,\ t_3=s_3s_5,\ t_4=s_4,\ t_5=s_5\rangle, \\
M_{e,i} &= \langle x,y\mid w^3=1,\ y^3=s_5^2,\ s_2^3=s_4^2s_5,\ s_3^3=s_5^2,\ t_3=s_3s_5w^{i-1},\ t_4=s_4,\ t_5=s_5\rangle.
\end{aligned}
\end{equation}
The \(p\)-class is given by
\(\mathrm{cl}_p(M_{e,i})=e+1\) respectively \(\mathrm{cl}_p(T_{e+1})=e+1\).
For \(e\ge 5\), the last non-trivial lower \(p\)-central is given by
\(P_{e}(M_{e,i})=\langle w\rangle\) respectively \(P_{e}(T_{e+1})=\langle w\rangle\),
whence all three groups share the common \(p\)-parent
\(M_{e,i}/P_{e}(M_{e,i})\simeq T_{e+1}/P_{e}(T_{e+1})\simeq T_e\).

Repeated recursive applications of the \(p\)-group generation algorithm
eventually produce a graph isomorphic to the tree in Figure
\ref{fig:SchurSigma5},
and thus confirm Formulas
\eqref{eqn:a1GroundD5},
\eqref{eqn:GroundD5M}
and
\eqref{eqn:GroundD5}.
All groups have \(\varrho\sim (2,2,3;3)\).
Cf.
\cite[Thm. 11--13]{Ma2021a}.
\end{proof}


Finally, we come to the ground state of type \(\mathrm{D}.6\).

\begin{theorem}
\label{thm:a1ChainGroundD6}
For each logarithmic exponent \(e\ge 5\),
a metabelian \(\mathrm{CF}\)-group \(T_e\)
with commutator quotient \(T_e/T_e^\prime\simeq (3^e,3)\hat{=}(e1)\),
punctured transfer kernel type \(\mathrm{a}.1\), \(\varkappa(T_e)\sim (000;0)\),
and rank distribution \(\varrho\sim (2,2,2;3)\)
is given by the \textbf{periodic sequence} of iterated \(p\)-descendants
\begin{equation}
\label{eqn:a1GroundD6}
T_e\simeq\mathrm{SmallGroup}(6561,85)-\#2;4(-\#1;1)^{e-5}.
\end{equation}
These groups have logarithmic order \(\mathrm{lo}(T_e)=5+e\) and
first abelian quotient invariants \(\alpha_1(T_e)\sim (e1,e1,e1;(e-1)22)\).
They form the infinite main trunk of a descendant tree
with finite double-twigs of depth two.
For each integer \(e\ge 5\), the doublet
\begin{equation}
\label{eqn:GroundD6M}
M_{e,i}\simeq\mathrm{SmallGroup}(6561,85)-\#2;4(-\#1;1)^{e-5}-\#1;i, \qquad i\in\lbrace 2,3\rbrace,
\end{equation}
respectively
\begin{equation}
\label{eqn:GroundD6}
G_{e,i}\simeq\mathrm{SmallGroup}(6561,85)-\#2;4(-\#1;1)^{e-5}-\#1;i-\#1;1, \qquad i\in\lbrace 2,3\rbrace,
\end{equation}
is the \textbf{unique pair} of \textbf{metabelian} \(\sigma\)-groups \(M_{e,i}\simeq G_{e,i}/G_{e,i}^{\prime\prime}\),
respectively \textbf{non-metabelian} Schur \(\sigma\)-groups \(G_{e,i}\), with
commutator quotient \((3^e,3)\hat{=}(e1)\),
punctured transfer kernel type \(\mathrm{D}.6\), \(\varkappa\sim (123;1)\),
\(\mathrm{lo}(M_{e,i})=6+e\), respectively \(\mathrm{lo}(G_{e,i})=7+e\), and
first abelian quotient invariants \(\alpha_1\sim ((e+1)1,(e+1)1,(e+1)1;e22)\).
See Figure
\ref{fig:SchurSigma5D6}.
\end{theorem}


\begin{proof}
Let \(s_2=\lbrack y,x\rbrack\) denote the main commutator,
\(s_3=\lbrack s_2,x\rbrack\), \(s_4=\lbrack s_3,x\rbrack\), \(s_5=\lbrack s_4,x\rbrack\)
and \(t_3=\lbrack s_2,y\rbrack\) higher commutators,
and \(w=x^{3^e}\) the last non-trivial power.
For each \(e\ge 5\), we have parametrized presentations
\begin{equation}
\label{eqn:Pres5D6}
\begin{aligned}
T_{e+1} &= \langle x,y\mid w^3=1,\ y^3=s_3^2s_4,\ s_2^3=s_4^2s_5,\ s_3^3=s_5^2,\ t_3=s_5\rangle, \\
M_{e,i} &= \langle x,y\mid w^3=1,\ y^3=s_3^2s_4,\ s_2^3=s_4^2s_5,\ s_3^3=s_5^2,\ t_3=s_5w^{i-1}\rangle, \qquad i\in\lbrace 2,3\rbrace.
\end{aligned}
\end{equation}
The \(p\)-class is given by
\(\mathrm{cl}_p(M_{e,i})=e+1\) respectively \(\mathrm{cl}_p(T_{e+1})=e+1\).
For \(e\ge 5\), the last non-trivial lower \(p\)-central is given by
\(P_{e}(M_{e,i})=\langle w\rangle\) respectively \(P_{e}(T_{e+1})=\langle w\rangle\),
whence all three groups share the common \(p\)-parent
\(M_{e,i}/P_{e}(M_{e,i})\simeq T_{e+1}/P_{e}(T_{e+1})\simeq T_e\).

Repeated recursive applications of the \(p\)-group generation algorithm
eventually produce Figure
\ref{fig:SchurSigma5D6},
and thus confirm Formulas
\eqref{eqn:a1GroundD6},
\eqref{eqn:GroundD6M}
and
\eqref{eqn:GroundD6}.
All groups have \(\varrho\sim (2,2,2;3)\).
\end{proof}


\begin{figure}[hb]
\caption{Schur \(\sigma\)-groups \(G\) with \(\varrho(G)\sim (2,2,3;3)\), \(G/G^\prime\simeq (3^e,3)\), \(2\le e\le 9\)}
\label{fig:SchurSigma7}

{\tiny

\setlength{\unitlength}{0.9cm}
\begin{picture}(14,18.5)(-9.5,-16.3)

\put(-11,2.5){\makebox(0,0)[cb]{order}}

\put(-11,2){\line(0,-1){17}}
\multiput(-11.1,2)(0,-1){18}{\line(1,0){0.2}}


\put(-10.8,2){\makebox(0,0)[lc]{\(3^2\)}}
\put(-10.8,1){\makebox(0,0)[lc]{\(3^3\)}}
\put(-10.8,0){\makebox(0,0)[lc]{\(3^4\)}}
\put(-10.8,-1){\makebox(0,0)[lc]{\(3^5\)}}
\put(-10.8,-2){\makebox(0,0)[lc]{\(3^6\)}}
\put(-10.8,-3){\makebox(0,0)[lc]{\(3^7\)}}
\put(-10.8,-4){\makebox(0,0)[lc]{\(3^8\)}}
\put(-10.8,-5){\makebox(0,0)[lc]{\(3^9\)}}
\put(-10.8,-6){\makebox(0,0)[lc]{\(3^{10}\)}}
\put(-10.8,-7){\makebox(0,0)[lc]{\(3^{11}\)}}
\put(-10.8,-8){\makebox(0,0)[lc]{\(3^{12}\)}}
\put(-10.8,-9){\makebox(0,0)[lc]{\(3^{13}\)}}
\put(-10.8,-10){\makebox(0,0)[lc]{\(3^{14}\)}}
\put(-10.8,-11){\makebox(0,0)[lc]{\(3^{15}\)}}
\put(-10.8,-12){\makebox(0,0)[lc]{\(3^{16}\)}}
\put(-10.8,-13){\makebox(0,0)[lc]{\(3^{17}\)}}
\put(-10.8,-14){\makebox(0,0)[lc]{\(3^{18}\)}}
\put(-10.8,-15){\makebox(0,0)[lc]{\(3^{19}\)}}


\put(-11,-15){\vector(0,-1){1.5}}

\put(-9,2){\circle{0.2}}
\put(-9,2){\circle*{0.1}}

\put(-9,1){\circle{0.2}}
\put(-9,-1){\circle{0.2}}
\put(-9,-2){\circle*{0.2}}
\put(-9,-3){\circle{0.2}}
\put(-9,-4){\circle{0.2}}
\put(-9,-5){\circle{0.2}}

\put(-7,0){\circle{0.2}}
\put(-7,-2){\circle{0.2}}
\put(-7,-3){\circle*{0.2}}
\put(-7,-4){\circle{0.2}}
\put(-7,-5){\circle{0.2}}
\put(-7,-6){\circle{0.2}}

\put(-5,-2){\circle{0.2}}
\put(-5,-4){\circle*{0.2}}
\put(-5,-5){\circle{0.2}}
\put(-5,-6){\circle{0.2}}
\put(-5,-7){\circle{0.2}}

\put(-3,-4){\circle{0.2}}
\put(-3,-6){\circle*{0.2}}
\put(-3,-7){\circle{0.2}}
\put(-3,-8){\circle{0.2}}

\put(-1,-6){\circle{0.2}}
\put(-1,-8){\circle*{0.2}}
\put(-1,-9){\circle{0.2}}

\put(1,-8){\circle{0.2}}
\put(1,-10){\circle{0.2}}

\put(3,-10){\circle{0.2}}
\put(3,-11){\circle{0.2}}

\put(4,-11){\circle{0.2}}
\put(4,-12){\circle{0.2}}

\put(5,-12){\circle{0.2}}
\put(5,-13){\circle{0.2}}

\put(-9.6,-4.1){\framebox(0.2,0.2){}}
\put(-9.6,-5.1){\framebox(0.2,0.2){}}
\put(-9.5,-5){\circle*{0.1}}
\put(-9.6,-7.1){\framebox(0.2,0.2){}}
\put(-9.5,-7){\circle{0.1}}

\put(-7.6,-5.1){\framebox(0.2,0.2){}}
\put(-7.6,-6.1){\framebox(0.2,0.2){}}
\put(-7.5,-6){\circle*{0.1}}
\put(-7.6,-8.1){\framebox(0.2,0.2){}}
\put(-7.5,-8){\circle{0.1}}

\put(-5.6,-6.1){\framebox(0.2,0.2){}}
\put(-5.6,-7.1){\framebox(0.2,0.2){}}
\put(-5.5,-7){\circle*{0.1}}
\put(-5.6,-9.1){\framebox(0.2,0.2){}}
\put(-5.5,-9){\circle{0.1}}

\put(-3.6,-8.1){\framebox(0.2,0.2){}}
\put(-3.5,-8){\circle*{0.1}}
\put(-3.6,-10.1){\framebox(0.2,0.2){}}
\put(-3.5,-10){\circle{0.1}}

\put(-1.6,-10.1){\framebox(0.2,0.2){}}
\put(-1.6,-11.1){\framebox(0.2,0.2){}}
\put(-1.5,-11){\circle{0.1}}

\put(0.9,-11.1){\framebox(0.2,0.2){}}
\put(0.9,-12.1){\framebox(0.2,0.2){}}
\put(1,-12){\circle{0.1}}

\put(2.9,-12.1){\framebox(0.2,0.2){}}
\put(2.9,-13.1){\framebox(0.2,0.2){}}
\put(3,-13){\circle{0.1}}

\put(3.9,-13.1){\framebox(0.2,0.2){}}
\put(3.9,-14.1){\framebox(0.2,0.2){}}
\put(4,-14){\circle{0.1}}

\put(4.9,-14.1){\framebox(0.2,0.2){}}
\put(4.9,-15.1){\framebox(0.2,0.2){}}
\put(5,-15){\circle{0.1}}


\put(-9,2){\line(0,-1){1}}
\put(-9,1){\line(0,-1){2}}
\put(-9,-1){\line(0,-1){1}}
\put(-9,-2){\line(0,-1){1}}
\put(-9,-3){\line(0,-1){1}}
\put(-9,-4){\line(0,-1){1}}
\put(-9,-2){\line(-1,-4){0.5}}
\put(-9.5,-4){\line(0,-1){1}}
\put(-9.5,-5){\line(0,-1){2}}

\put(-9,2){\line(1,-1){2}}
\put(-7,0){\line(0,-1){2}}
\put(-7,-2){\line(0,-1){1}}
\put(-7,-3){\line(0,-1){1}}
\put(-7,-3){\line(-1,-4){0.5}}
\put(-7,-4){\line(0,-1){1}}
\put(-7,-5){\line(0,-1){1}}
\put(-7.5,-5){\line(0,-1){1}}
\put(-7.5,-6){\line(0,-1){2}}

\put(-7,0){\line(1,-1){2}}
\put(-5,-2){\line(0,-1){2}}
\put(-5,-4){\line(0,-1){1}}
\put(-5,-5){\line(0,-1){1}}
\put(-5,-6){\line(0,-1){1}}
\put(-5,-4){\line(-1,-4){0.5}}
\put(-5.5,-6){\line(0,-1){1}}
\put(-5.5,-7){\line(0,-1){2}}

\put(-5,-2){\line(1,-1){2}}
\put(-3,-4){\line(0,-1){2}}
\put(-3,-6){\line(-1,-4){0.5}}
\put(-3,-6){\line(0,-1){1}}
\put(-3,-7){\line(0,-1){1}}
\put(-3.5,-8){\line(0,-1){2}}

\put(-3,-4){\line(1,-1){2}}
\put(-1,-6){\line(0,-1){2}}
\put(-1,-8){\line(0,-1){1}}
\put(-1,-8){\line(-1,-4){0.5}}
\put(-1.5,-10){\line(0,-1){1}}

\put(-1,-6){\line(1,-1){2}}
\put(1,-8){\line(0,-1){2}}
\put(1,-10){\line(0,-1){1}}
\put(1,-11){\line(0,-1){1}}

\put(1,-8){\line(1,-1){2}}
\put(3,-10){\line(0,-1){1}}
\put(3,-11){\line(0,-1){1}}
\put(3,-12){\line(0,-1){1}}

\put(3,-10){\line(1,-1){1}}
\put(4,-11){\line(0,-1){1}}
\put(4,-12){\line(0,-1){1}}
\put(4,-13){\line(0,-1){1}}

\put(4,-11){\line(1,-1){1}}
\put(5,-12){\line(0,-1){1}}
\put(5,-13){\line(0,-1){1}}
\put(5,-14){\line(0,-1){1}}

\put(5,-12){\line(1,-1){1}}


\put(-8.8,2){\makebox(0,0)[lc]{\(\langle 2\rangle\)}}

\put(-8.8,1){\makebox(0,0)[lc]{\(\langle 3\rangle\)}}
\put(-8.8,-1){\makebox(0,0)[lc]{\(\langle 6\rangle\)}}
\put(-8.8,-2){\makebox(0,0)[lc]{\(\langle 49\rangle\)}}
\put(-8.8,-3){\makebox(0,0)[lc]{\(1;2\)}}
\put(-8.8,-4){\makebox(0,0)[lc]{\(1;1\)}}
\put(-8.8,-5){\makebox(0,0)[lc]{\(1;i\)}}

\put(-6.8,0){\makebox(0,0)[lc]{\(\langle 3\rangle\)}}
\put(-6.8,-2){\makebox(0,0)[lc]{\(\langle 13\rangle\)}}
\put(-6.8,-3){\makebox(0,0)[lc]{\(\langle 168\rangle\)}}
\put(-6.8,-4){\makebox(0,0)[lc]{\(1;7\)}}
\put(-6.8,-5){\makebox(0,0)[lc]{\(1;4\)}}
\put(-6.8,-6){\makebox(0,0)[lc]{\(1;i\)}}

\put(-4.8,-2){\makebox(0,0)[lc]{\(\langle 7\rangle\)}}
\put(-4.8,-4){\makebox(0,0)[lc]{\(\langle 98\rangle\)}}
\put(-4.8,-5){\makebox(0,0)[lc]{\(1;3\)}}
\put(-4.8,-6){\makebox(0,0)[lc]{\(1;1\)}}
\put(-4.8,-7){\makebox(0,0)[lc]{\(1;i\)}}

\put(-2.8,-4){\makebox(0,0)[lc]{\(\langle 93\rangle\)}}
\put(-2.8,-6){\makebox(0,0)[lc]{\(2;6\)}}
\put(-2.8,-7){\makebox(0,0)[lc]{\(1;2\)}}
\put(-2.8,-8){\makebox(0,0)[lc]{\(1;i\)}}

\put(-0.8,-6){\makebox(0,0)[lc]{\(2;1\)}}
\put(-0.8,-8){\makebox(0,0)[lc]{\(2;6\)}}
\put(-0.8,-9){\makebox(0,0)[lc]{\(1;i+1\)}}

\put(1.2,-6.5){\vector(-1,-1){1.5}}
\put(1.3,-6.5){\makebox(0,0)[lc]{last semi-metabelian bifurcation}}

\put(1.2,-8){\makebox(0,0)[lc]{\(2;1\)}}
\put(1.2,-10){\makebox(0,0)[lc]{\(2;i\)}}

\put(3.2,-10){\makebox(0,0)[lc]{\(2;2\), periodic root \(T_7\)}}
\put(3.8,-10.3){\makebox(0,0)[lc]{in Formula \eqref{eqn:b16ExcitedD10}}}
\put(3.2,-11){\makebox(0,0)[lc]{\(1;i\)}}

\put(4.2,-11){\makebox(0,0)[lc]{\(1;1\)}}
\put(4.2,-12){\makebox(0,0)[lc]{\(1;i\)}}

\put(5.2,-12){\makebox(0,0)[lc]{\(1;1\)}}
\put(5.2,-13){\makebox(0,0)[lc]{\(1;i\)}}

\put(-9.7,-4){\makebox(0,0)[rc]{\(2;1\)}}
\put(-9.7,-5){\makebox(0,0)[rc]{\(1;1\)}}
\put(-9.7,-7){\makebox(0,0)[rc]{\(2;i\)}}
\put(-9.5,-7.5){\makebox(0,0)[cc]{\(i=5,6\)}}
\put(-9.5,-8){\makebox(0,0)[lc]{\(\mathrm{E}.14\)}}
\put(-9.5,-8.5){\makebox(0,0)[lc]{\((3,3)\)}}

\put(-7.7,-5){\makebox(0,0)[rc]{\(2;7\)}}
\put(-7.7,-6){\makebox(0,0)[rc]{\(1;4\)}}
\put(-7.7,-8){\makebox(0,0)[rc]{\(2;i\)}}
\put(-7.5,-8.5){\makebox(0,0)[cc]{\(i=5,6\)}}
\put(-7.5,-9){\makebox(0,0)[lc]{\(\mathrm{D}.10\)}}
\put(-7.5,-9.5){\makebox(0,0)[lc]{\((9,3)\)}}

\put(-5.7,-6){\makebox(0,0)[rc]{\(2;1\)}}
\put(-5.7,-7){\makebox(0,0)[rc]{\(1;1\)}}
\put(-5.7,-9){\makebox(0,0)[rc]{\(2;i\)}}
\put(-5.5,-9.5){\makebox(0,0)[cc]{\(i=2,3\)}}
\put(-5.5,-10){\makebox(0,0)[lc]{\(\mathrm{D}.10\)}}
\put(-5.5,-10.5){\makebox(0,0)[lc]{\((27,3)\)}}

\put(-3.7,-8){\makebox(0,0)[rc]{\(2;1\)}}
\put(-3.7,-10){\makebox(0,0)[rc]{\(2;i\)}}
\put(-3.5,-10.5){\makebox(0,0)[cc]{\(i=2,3\)}}
\put(-3.5,-11){\makebox(0,0)[lc]{\(\mathrm{D}.10\)}}
\put(-3.5,-11.5){\makebox(0,0)[lc]{\((81,3)\)}}

\put(-2.1,-13.5){\vector(-1,4){1.3}}
\put(-2.3,-13.6){\makebox(0,0)[ct]{last non-metabelian bifurcation}}

\put(-1.7,-10){\makebox(0,0)[rc]{\(2;i\)}}
\put(-1.7,-11){\makebox(0,0)[rc]{\(1;1\)}}
\put(-1.5,-11.5){\makebox(0,0)[cc]{\(i=2,3\)}}
\put(-1.5,-12){\makebox(0,0)[lc]{\(\mathrm{D}.10\)}}
\put(-1.5,-12.5){\makebox(0,0)[lc]{\((243,3)\)}}

\put(1.2,-11){\makebox(0,0)[lc]{\(1;1\)}}
\put(1.2,-12){\makebox(0,0)[lc]{\(1;1\)}}
\put(0.5,-12.5){\makebox(0,0)[lc]{\(i=7,8\)}}
\put(0.5,-13){\makebox(0,0)[lc]{\(\mathrm{D}.10\)}}
\put(0.5,-13.5){\makebox(0,0)[lc]{\((729,3)\)}}

\put(3.2,-12){\makebox(0,0)[lc]{\(1;1\)}}
\put(3.2,-13){\makebox(0,0)[lc]{\(1;1\)}}
\put(2.5,-13.5){\makebox(0,0)[lc]{\(i=2,3\)}}
\put(2.5,-14){\makebox(0,0)[lc]{\(\mathrm{D}.10\)}}
\put(2.5,-14.5){\makebox(0,0)[lc]{\((2187,3)\)}}

\put(4.2,-13){\makebox(0,0)[lc]{\(1;1\)}}
\put(4.2,-14){\makebox(0,0)[lc]{\(1;1\)}}
\put(3.8,-14.5){\makebox(0,0)[lc]{\(i=2,3\)}}
\put(3.5,-15){\makebox(0,0)[lc]{\(\mathrm{D}.10\)}}
\put(3.5,-15.5){\makebox(0,0)[lc]{\((6561,3)\)}}

\put(5.2,-14){\makebox(0,0)[lc]{\(1;1\)}}
\put(5.2,-15){\makebox(0,0)[lc]{\(1;1\)}}
\put(4.8,-15.5){\makebox(0,0)[lc]{\(i=2,3\)}}
\put(4.5,-16){\makebox(0,0)[lc]{\(\mathrm{D}.10\)}}
\put(4.5,-16.5){\makebox(0,0)[lc]{\((19683,3)\)}}

\put(4.8,-9.5){\makebox(0,0)[cc]{main trunk of type \(\mathrm{b}.16\)}}
\put(6.2,-13){\makebox(0,0)[lc]{etc.}}

\put(-4,1.4){\makebox(0,0)[lc]{Legend:}}

\put(-2.5,1.4){\circle{0.2}}
\put(-2.5,1.4){\circle*{0.1}}
\put(-2,1.4){\makebox(0,0)[lc]{\(\ldots\) abelian}}

\put(-2.5,1){\circle{0.2}}
\put(-2,1){\makebox(0,0)[lc]{\(\ldots\) metabelian}}

\put(-2.5,0.6){\circle*{0.2}}
\put(-2,0.6){\makebox(0,0)[lc]{\(\ldots\) metabelian with bifurcation}}

\put(-2.6,0.1){\framebox(0.2,0.2){}}
\put(-2,0.2){\makebox(0,0)[lc]{\(\ldots\) non-metabelian}}

\put(-2.6,-0.3){\framebox(0.2,0.2){}}
\put(-2.5,-0.2){\circle*{0.1}}
\put(-2,-0.2){\makebox(0,0)[lc]{\(\ldots\) non-metabelian with bifurcation}}

\put(-2.6,-0.7){\framebox(0.2,0.2){}}
\put(-2.5,-0.6){\circle{0.1}}
\put(-2,-0.6){\makebox(0,0)[lc]{\(\ldots\) non-metabelian Schur \(\sigma\)}}

\end{picture}

}

\end{figure}


\section{Moderate Schur \(\sigma\)-groups \(G\) with \(G/G^\prime\simeq (3^e,3)\), \(e\ge 7\), excited state}
\label{s:NonMetabelian7}


We begin with a thorough investigation of the first excited state of type \(\mathrm{D}.10\).

\begin{theorem}
\label{thm:b16ChainExcitedD10}
For each logarithmic exponent \(e\ge 7\),
a metabelian \(\mathrm{CF}\)-group \(T_e\)
with commutator quotient \(T_e/T_e^\prime\simeq (3^e,3)\hat{=}(e1)\),
punctured transfer kernel type \(\mathrm{b}.16\), \(\varkappa(T_e)\sim (004;0)\),
and rank distribution \(\varrho\sim (2,2,3;3)\)
is given by the \textbf{periodic sequence} of iterated \(p\)-descendants
\begin{equation}
\label{eqn:b16ExcitedD10}
T_e\simeq\mathrm{SmallGroup}(6561,93)(-\#2;1)^2-\#2;2(-\#1;1)^{e-7}.
\end{equation}
These groups have logarithmic order \(\mathrm{lo}(T_e)=7+e\) and
first abelian quotient invariants \(\alpha_1(T_e)\sim (e1,e1,e33;(e-1)11)\).
They form the infinite main trunk of a descendant tree
with finite double-twigs of depth three.
For each integer \(e\ge 7\), the doublet
\begin{equation}
\label{eqn:ExcitedD10M}
M_{e,i}\simeq\mathrm{SmallGroup}(6561,93)(-\#2;1)^2-\#2;2(-\#1;1)^{e-7}-\#1;i, \qquad i\in\lbrace 2,3\rbrace,
\end{equation}
respectively
\begin{equation}
\label{eqn:ExcitedD10}
G_{e,i}\simeq\mathrm{SmallGroup}(6561,93)(-\#2;1)^2-\#2;2(-\#1;1)^{e-7}-\#1;i(-\#1;1)^2, \qquad i\in\lbrace 2,3\rbrace,
\end{equation}
is the \textbf{unique pair} of \textbf{metabelian} \(\sigma\)-groups \(M_{e,i}\simeq G_{e,i}/G_{e,i}^{\prime\prime}\),
respectively \textbf{non-metabelian} Schur \(\sigma\)-groups \(G_{e,i}\), with
commutator quotient \((3^e,3)\hat{=}(e1)\),
punctured transfer kernel type \(\mathrm{D}.10\), \(\varkappa\sim (114;3)\),
\(\mathrm{lo}(M_{e,i})=8+e\), respectively \(\mathrm{lo}(G_{e,i})=10+e\), and
first abelian quotient invariants \(\alpha_1\sim ((e+1)1,(e+1)1,e33;e11)\).
See Figure
\ref{fig:SchurSigma7}.
\end{theorem}


\begin{proof}
Let \(s_2=\lbrack y,x\rbrack\) denote the main commutator,
\(\forall_{j=3}^7\) \(s_j=\lbrack s_{j-1},x\rbrack\) 
and \(t_j=\lbrack s_{j-1},y\rbrack\) higher commutators,
and \(w=x^{3^e}\) the last non-trivial power.
For each \(e\ge 7\), we have parametrized presentations with \(i\in\lbrace 2,3\rbrace\)
\begin{equation}
\label{eqn:Pres7D10}
\begin{aligned}
T_{e+1} &= \langle x,y\mid w^3=1,\ y^3=s_7,\ \forall_{j=2}^4\,s_j^3=s_{j+2}^2s_{j+3},\ s_5^3=s_7^2,\ \forall_{j=3}^7\,t_j=s_j\rangle, \\
M_{e,i} &= \langle x,y\mid w^3=1,\ y^3=s_7,\ \forall_{j=2}^4\,s_j^3=s_{j+2}^2s_{j+3},\ s_5^3=s_7^2,\ t_3=s_3w^{i-1},\ \forall_{j=4}^7\,t_j=s_j\rangle.
\end{aligned}
\end{equation}
The \(p\)-class is given by
\(\mathrm{cl}_p(M_{e,i})=e+1\) respectively \(\mathrm{cl}_p(T_{e+1})=e+1\).
For \(e\ge 7\), the last non-trivial lower \(p\)-central is given by
\(P_{e}(M_{e,i})=\langle w\rangle\) respectively \(P_{e}(T_{e+1})=\langle w\rangle\),
whence all three groups share the common \(p\)-parent
\(M_{e,i}/P_{e}(M_{e,i})\simeq T_{e+1}/P_{e}(T_{e+1})\simeq T_e\).

Repeated recursive applications of the \(p\)-group generation algorithm
eventually produce Figure
\ref{fig:SchurSigma7},
and thus confirm Formulas
\eqref{eqn:b16ExcitedD10},
\eqref{eqn:ExcitedD10M}
and
\eqref{eqn:ExcitedD10}.
All groups have \(\varrho\sim (2,2,3;3)\).
Cf.
\cite[Thm. 1--3]{Ma2021c}.
\end{proof}


We continue with the first excited state of type \(\mathrm{C}.4\).

\begin{theorem}
\label{thm:a1ChainExcitedC4}
For each logarithmic exponent \(e\ge 7\),
a metabelian \(\mathrm{CF}\)-group \(T_e\)
with commutator quotient \(T_e/T_e^\prime\simeq (3^e,3)\hat{=}(e1)\),
punctured transfer kernel type \(\mathrm{a}.1\), \(\varkappa(T_e)\sim (000;0)\),
and rank distribution \(\varrho\sim (2,2,3;3)\)
is given by the \textbf{periodic sequence} of iterated \(p\)-descendants
\begin{equation}
\label{eqn:a1ExcitedC4}
T_e\simeq\mathrm{SmallGroup}(6561,93)(-\#2;1)^2-\#2;5(-\#1;1)^{e-7}.
\end{equation}
These groups have logarithmic order \(\mathrm{lo}(T_e)=7+e\) and
first abelian quotient invariants \(\alpha_1(T_e)\sim (e1,e1,e32;(e-1)11)\).
They form the infinite main trunk of a descendant tree
with finite double-twigs of depth three.
For each integer \(e\ge 7\), the doublet
\begin{equation}
\label{eqn:ExcitedC4M}
M_{e,i}\simeq\mathrm{SmallGroup}(6561,93)(-\#2;1)^2-\#2;5(-\#1;1)^{e-7}-\#1;i, \qquad i\in\lbrace 2,3\rbrace,
\end{equation}
respectively
\begin{equation}
\label{eqn:ExcitedC4}
G_{e,i}\simeq\mathrm{SmallGroup}(6561,93)(-\#2;1)^2-\#2;5(-\#1;1)^{e-7}-\#1;i(-\#1;1)^2, \qquad i\in\lbrace 2,3\rbrace,
\end{equation}
is the \textbf{unique pair} of \textbf{metabelian} \(\sigma\)-groups \(M_{e,i}\simeq G_{e,i}/G_{e,i}^{\prime\prime}\),
respectively \textbf{non-metabelian} Schur \(\sigma\)-groups \(G_{e,i}\), with
commutator quotient \((3^e,3)\hat{=}(e1)\),
punctured transfer kernel type \(\mathrm{C}.4\), \(\varkappa\sim (113;3)\),
\(\mathrm{lo}(M_{e,i})=8+e\), respectively \(\mathrm{lo}(G_{e,i})=10+e\), and
first abelian quotient invariants \(\alpha_1\sim ((e+1)1,(e+1)1,(e+1)32;e11)\).
\end{theorem}


\begin{proof}
Let \(s_2=\lbrack y,x\rbrack\) denote the main commutator,
\(\forall_{j=3}^7\) \(s_j=\lbrack s_{j-1},x\rbrack\)
and \(t_j=\lbrack s_{j-1},y\rbrack\) higher commutators,
and \(w=x^{3^e}\) the last non-trivial power.
For each \(e\ge 7\), we have parametrized presentations with \(i\in\lbrace 2,3\rbrace\)
\begin{equation}
\label{eqn:Pres7C4}
\begin{aligned}
T_{e+1} &= \langle x,y\mid w^3=1,\ y^3=s_7,\ \forall_{j=2}^4\,s_j^3=s_{j+2}^2s_{j+3},s_5^3=s_7^2,\ t_3=s_3s_7,\ \forall_{j=4}^7\,t_j=s_j\rangle, \\
M_{e,i} &= \langle x,y\mid w^3=1,\ y^3=s_7,\ \forall_{j=2}^4\,s_j^3=s_{j+2}^2s_{j+3},s_5^3=s_7^2,\ t_3=s_3s_7w^{i-1},\ \forall_{j=4}^7\,t_j=s_j\rangle.
\end{aligned}
\end{equation}
The \(p\)-class is given by
\(\mathrm{cl}_p(M_{e,i})=e+1\) respectively \(\mathrm{cl}_p(T_{e+1})=e+1\).
For \(e\ge 7\), the last non-trivial lower \(p\)-central is given by
\(P_{e}(M_{e,i})=\langle w\rangle\) respectively \(P_{e}(T_{e+1})=\langle w\rangle\),
whence all three groups share the common \(p\)-parent
\(M_{e,i}/P_{e}(M_{e,i})\simeq T_{e+1}/P_{e}(T_{e+1})\simeq T_e\).

Repeated recursive applications of the \(p\)-group generation algorithm
eventually produce a graph isomorphic to the tree in Figure
\ref{fig:SchurSigma7},
and thus confirm Formulas
\eqref{eqn:a1ExcitedC4},
\eqref{eqn:ExcitedC4M}
and
\eqref{eqn:ExcitedC4}.
All groups have \(\varrho\sim (2,2,3;3)\).
Cf.
\cite[Thm. 1--3]{Ma2021c}.
\end{proof}


Next, we look at the first excited state of type \(\mathrm{D}.5\).

\begin{theorem}
\label{thm:a1ChainExcitedD5}
For each logarithmic exponent \(e\ge 7\),
a metabelian \(\mathrm{CF}\)-group \(T_e\)
with commutator quotient \(T_e/T_e^\prime\simeq (3^e,3)\hat{=}(e1)\),
punctured transfer kernel type \(\mathrm{a}.1\), \(\varkappa(T_e)\sim (000;0)\),
and rank distribution \(\varrho\sim (2,2,3;3)\)
is given by the \textbf{periodic sequence} of iterated \(p\)-descendants
\begin{equation}
\label{eqn:a1ExcitedD5}
T_e\simeq\mathrm{SmallGroup}(6561,93)(-\#2;1)^2-\#2;4(-\#1;1)^{e-7}.
\end{equation}
These groups have logarithmic order \(\mathrm{lo}(T_e)=7+e\) and
first abelian quotient invariants \(\alpha_1(T_e)\sim (e1,e1,e32;(e-1)11)\).
They form the infinite main trunk of a descendant tree
with finite double-twigs of depth three.
For each integer \(e\ge 7\), the doublet
\begin{equation}
\label{eqn:ExcitedD5M}
M_{e,i}\simeq\mathrm{SmallGroup}(6561,93)(-\#2;1)^2-\#2;4(-\#1;1)^{e-7}-\#1;i, \qquad i\in\lbrace 2,3\rbrace,
\end{equation}
respectively
\begin{equation}
\label{eqn:ExcitedD5}
G_{e,i}\simeq\mathrm{SmallGroup}(6561,93)(-\#2;1)^2-\#2;4(-\#1;1)^{e-7}-\#1;i(-\#1;1)^2, \qquad i\in\lbrace 2,3\rbrace,
\end{equation}
is the \textbf{unique pair} of \textbf{metabelian} \(\sigma\)-groups \(M_{e,i}\simeq G_{e,i}/G_{e,i}^{\prime\prime}\),
respectively \textbf{non-metabelian} Schur \(\sigma\)-groups \(G_{e,i}\), with
commutator quotient \((3^e,3)\hat{=}(e1)\),
punctured transfer kernel type \(\mathrm{D}.5\), \(\varkappa\sim (112;3)\),
\(\mathrm{lo}(M_{e,i})=8+e\), respectively \(\mathrm{lo}(G_{e,i})=10+e\), and
first abelian quotient invariants \(\alpha_1\sim ((e+1)1,(e+1)1,(e+1)32;e11)\).
\end{theorem}


\begin{proof}
Let \(s_2=\lbrack y,x\rbrack\) denote the main commutator,
\(\forall_{j=3}^7\) \(s_j=\lbrack s_{j-1},x\rbrack\)
and \(t_j=\lbrack s_{j-1},y\rbrack\) higher commutators,
and \(w=x^{3^e}\) the last non-trivial power.
For each \(e\ge 7\), we have parametrized presentations with \(i\in\lbrace 2,3\rbrace\)
\begin{equation}
\label{eqn:Pres7D5}
\begin{aligned}
T_{e+1} &= \langle x,y\mid w^3=1,\ y^3=s_7^2,\ \forall_{j=2}^4\,s_j^3=s_{j+2}^2s_{j+3},s_5^3=s_7^2,\ t_3=s_3s_7,\ \forall_{j=4}^7\,t_j=s_j\rangle, \\
M_{e,i} &= \langle x,y\mid w^3=1,\ y^3=s_7^2,\ \forall_{j=2}^4\,s_j^3=s_{j+2}^2s_{j+3},s_5^3=s_7^2,\ t_3=s_3s_7w^{i-1},\ \forall_{j=4}^7\,t_j=s_j\rangle.
\end{aligned}
\end{equation}
The \(p\)-class is given by
\(\mathrm{cl}_p(M_{e,i})=e+1\) respectively \(\mathrm{cl}_p(T_{e+1})=e+1\).
For \(e\ge 7\), the last non-trivial lower \(p\)-central is given by
\(P_{e}(M_{e,i})=\langle w\rangle\) respectively \(P_{e}(T_{e+1})=\langle w\rangle\),
whence all three groups share the common \(p\)-parent
\(M_{e,i}/P_{e}(M_{e,i})\simeq T_{e+1}/P_{e}(T_{e+1})\simeq T_e\).

Repeated recursive applications of the \(p\)-group generation algorithm
eventually produce a graph isomorphic to the tree in Figure
\ref{fig:SchurSigma7},
and thus confirm Formulas
\eqref{eqn:a1ExcitedD5},
\eqref{eqn:ExcitedD5M}
and
\eqref{eqn:ExcitedD5}.
All groups have \(\varrho\sim (2,2,3;3)\).
Cf.
\cite[Thm. 1--3]{Ma2021c}.
\end{proof}


\begin{figure}[hb]
\caption{Schur \(\sigma\)-groups \(G\) with \(\varrho(G)\sim (2,2,2;3)\), \(G/G^\prime\simeq (3^e,3)\), \(2\le e\le 9\)}
\label{fig:SchurSigma7D6}

{\tiny

\setlength{\unitlength}{0.9cm}
\begin{picture}(14,18.5)(-9.5,-16.3)

\put(-11,2.5){\makebox(0,0)[cb]{order}}

\put(-11,2){\line(0,-1){17}}
\multiput(-11.1,2)(0,-1){18}{\line(1,0){0.2}}


\put(-10.8,2){\makebox(0,0)[lc]{\(3^2\)}}
\put(-10.8,1){\makebox(0,0)[lc]{\(3^3\)}}
\put(-10.8,0){\makebox(0,0)[lc]{\(3^4\)}}
\put(-10.8,-1){\makebox(0,0)[lc]{\(3^5\)}}
\put(-10.8,-2){\makebox(0,0)[lc]{\(3^6\)}}
\put(-10.8,-3){\makebox(0,0)[lc]{\(3^7\)}}
\put(-10.8,-4){\makebox(0,0)[lc]{\(3^8\)}}
\put(-10.8,-5){\makebox(0,0)[lc]{\(3^9\)}}
\put(-10.8,-6){\makebox(0,0)[lc]{\(3^{10}\)}}
\put(-10.8,-7){\makebox(0,0)[lc]{\(3^{11}\)}}
\put(-10.8,-8){\makebox(0,0)[lc]{\(3^{12}\)}}
\put(-10.8,-9){\makebox(0,0)[lc]{\(3^{13}\)}}
\put(-10.8,-10){\makebox(0,0)[lc]{\(3^{14}\)}}
\put(-10.8,-11){\makebox(0,0)[lc]{\(3^{15}\)}}
\put(-10.8,-12){\makebox(0,0)[lc]{\(3^{16}\)}}
\put(-10.8,-13){\makebox(0,0)[lc]{\(3^{17}\)}}
\put(-10.8,-14){\makebox(0,0)[lc]{\(3^{18}\)}}
\put(-10.8,-15){\makebox(0,0)[lc]{\(3^{19}\)}}


\put(-11,-15){\vector(0,-1){1.5}}

\put(-9,2){\circle{0.2}}
\put(-9,2){\circle*{0.1}}

\put(-9,1){\circle{0.2}}
\put(-9,-1){\circle{0.2}}
\put(-9,-2){\circle*{0.2}}
\put(-9,-3){\circle{0.2}}
\put(-9,-4){\circle{0.2}}
\put(-9,-5){\circle{0.2}}

\put(-7,0){\circle{0.2}}
\put(-7,-2){\circle{0.2}}
\put(-7,-3){\circle*{0.2}}
\put(-7,-4){\circle{0.2}}
\put(-7,-5){\circle{0.2}}
\put(-7,-6){\circle{0.2}}

\put(-5,-2){\circle{0.2}}
\put(-5,-4){\circle*{0.2}}
\put(-5,-5){\circle{0.2}}
\put(-5,-6){\circle{0.2}}
\put(-5,-7){\circle{0.2}}

\put(-3,-4){\circle{0.2}}
\put(-3,-6){\circle*{0.2}}
\put(-3,-7){\circle{0.2}}
\put(-3,-8){\circle{0.2}}

\put(-1,-6){\circle{0.2}}
\put(-1,-8){\circle*{0.2}}
\put(-1,-9){\circle{0.2}}

\put(1,-8){\circle{0.2}}
\put(1,-10){\circle{0.2}}

\put(3,-10){\circle{0.2}}
\put(3,-11){\circle{0.2}}

\put(4,-11){\circle{0.2}}
\put(4,-12){\circle{0.2}}

\put(5,-12){\circle{0.2}}
\put(5,-13){\circle{0.2}}

\put(-9.6,-4.1){\framebox(0.2,0.2){}}
\put(-9.6,-5.1){\framebox(0.2,0.2){}}
\put(-9.5,-5){\circle*{0.1}}
\put(-9.6,-7.1){\framebox(0.2,0.2){}}
\put(-9.5,-7){\circle{0.1}}

\put(-7.6,-5.1){\framebox(0.2,0.2){}}
\put(-7.6,-6.1){\framebox(0.2,0.2){}}
\put(-7.5,-6){\circle*{0.1}}
\put(-7.6,-8.1){\framebox(0.2,0.2){}}
\put(-7.5,-8){\circle{0.1}}

\put(-5.6,-6.1){\framebox(0.2,0.2){}}
\put(-5.6,-7.1){\framebox(0.2,0.2){}}
\put(-5.5,-7){\circle*{0.1}}
\put(-5.6,-9.1){\framebox(0.2,0.2){}}
\put(-5.5,-9){\circle{0.1}}

\put(-3.6,-8.1){\framebox(0.2,0.2){}}
\put(-3.5,-8){\circle*{0.1}}
\put(-3.6,-10.1){\framebox(0.2,0.2){}}
\put(-3.5,-10){\circle{0.1}}

\put(-1.6,-10.1){\framebox(0.2,0.2){}}
\put(-1.6,-11.1){\framebox(0.2,0.2){}}
\put(-1.5,-11){\circle{0.1}}

\put(0.9,-11.1){\framebox(0.2,0.2){}}
\put(0.9,-12.1){\framebox(0.2,0.2){}}
\put(1,-12){\circle{0.1}}

\put(2.9,-12.1){\framebox(0.2,0.2){}}
\put(2.9,-13.1){\framebox(0.2,0.2){}}
\put(3,-13){\circle{0.1}}

\put(3.9,-13.1){\framebox(0.2,0.2){}}
\put(3.9,-14.1){\framebox(0.2,0.2){}}
\put(4,-14){\circle{0.1}}

\put(4.9,-14.1){\framebox(0.2,0.2){}}
\put(4.9,-15.1){\framebox(0.2,0.2){}}
\put(5,-15){\circle{0.1}}


\put(-9,2){\line(0,-1){1}}
\put(-9,1){\line(0,-1){2}}
\put(-9,-1){\line(0,-1){1}}
\put(-9,-2){\line(0,-1){1}}
\put(-9,-3){\line(0,-1){1}}
\put(-9,-4){\line(0,-1){1}}
\put(-9,-2){\line(-1,-4){0.5}}
\put(-9.5,-4){\line(0,-1){1}}
\put(-9.5,-5){\line(0,-1){2}}

\put(-9,2){\line(1,-1){2}}
\put(-7,0){\line(0,-1){2}}
\put(-7,-2){\line(0,-1){1}}
\put(-7,-3){\line(0,-1){1}}
\put(-7,-3){\line(-1,-4){0.5}}
\put(-7,-4){\line(0,-1){1}}
\put(-7,-5){\line(0,-1){1}}
\put(-7.5,-5){\line(0,-1){1}}
\put(-7.5,-6){\line(0,-1){2}}

\put(-7,0){\line(1,-1){2}}
\put(-5,-2){\line(0,-1){2}}
\put(-5,-4){\line(0,-1){1}}
\put(-5,-5){\line(0,-1){1}}
\put(-5,-6){\line(0,-1){1}}
\put(-5,-4){\line(-1,-4){0.5}}
\put(-5.5,-6){\line(0,-1){1}}
\put(-5.5,-7){\line(0,-1){2}}

\put(-5,-2){\line(1,-1){2}}
\put(-3,-4){\line(0,-1){2}}
\put(-3,-6){\line(-1,-4){0.5}}
\put(-3,-6){\line(0,-1){1}}
\put(-3,-7){\line(0,-1){1}}
\put(-3.5,-8){\line(0,-1){2}}

\put(-3,-4){\line(1,-1){2}}
\put(-1,-6){\line(0,-1){2}}
\put(-1,-8){\line(0,-1){1}}
\put(-1,-8){\line(-1,-4){0.5}}
\put(-1.5,-10){\line(0,-1){1}}

\put(-1,-6){\line(1,-1){2}}
\put(1,-8){\line(0,-1){2}}
\put(1,-10){\line(0,-1){1}}
\put(1,-11){\line(0,-1){1}}

\put(1,-8){\line(1,-1){2}}
\put(3,-10){\line(0,-1){1}}
\put(3,-11){\line(0,-1){1}}
\put(3,-12){\line(0,-1){1}}

\put(3,-10){\line(1,-1){1}}
\put(4,-11){\line(0,-1){1}}
\put(4,-12){\line(0,-1){1}}
\put(4,-13){\line(0,-1){1}}

\put(4,-11){\line(1,-1){1}}
\put(5,-12){\line(0,-1){1}}
\put(5,-13){\line(0,-1){1}}
\put(5,-14){\line(0,-1){1}}

\put(5,-12){\line(1,-1){1}}


\put(-8.8,2){\makebox(0,0)[lc]{\(\langle 2\rangle\)}}

\put(-8.8,1){\makebox(0,0)[lc]{\(\langle 3\rangle\)}}
\put(-8.8,-1){\makebox(0,0)[lc]{\(\langle 8\rangle\)}}
\put(-8.8,-2){\makebox(0,0)[lc]{\(\langle 54\rangle\)}}
\put(-8.8,-3){\makebox(0,0)[lc]{\(1;3\)}}
\put(-8.8,-4){\makebox(0,0)[lc]{\(1;1\)}}
\put(-8.8,-5){\makebox(0,0)[lc]{\(1;i\)}}

\put(-6.8,0){\makebox(0,0)[lc]{\(\langle 3\rangle\)}}
\put(-6.8,-2){\makebox(0,0)[lc]{\(\langle 18\vert 21\rangle\)}}
\put(-6.8,-3){\makebox(0,0)[lc]{\(\langle 181\vert 191\rangle\)}}
\put(-6.8,-4){\makebox(0,0)[lc]{\(1;1\)}}
\put(-6.8,-5){\makebox(0,0)[lc]{\(1;1\)}}
\put(-6.8,-6){\makebox(0,0)[lc]{\(1;4\)}}

\put(-4.8,-2){\makebox(0,0)[lc]{\(\langle 6\rangle\)}}
\put(-4.8,-4){\makebox(0,0)[lc]{\(\langle 88\vert 91\rangle\)}}
\put(-4.8,-5){\makebox(0,0)[lc]{\(1;1\)}}
\put(-4.8,-6){\makebox(0,0)[lc]{\(1;1\)}}
\put(-4.8,-7){\makebox(0,0)[lc]{\(1;4\)}}

\put(-2.8,-4){\makebox(0,0)[lc]{\(\langle 85\rangle\)}}
\put(-2.8,-6){\makebox(0,0)[lc]{\(2;5\vert 9\)}}
\put(-2.8,-7){\makebox(0,0)[lc]{\(1;2\)}}
\put(-2.8,-8){\makebox(0,0)[lc]{\(1;4\)}}

\put(-0.8,-6){\makebox(0,0)[lc]{\(2;1\)}}
\put(-0.8,-8){\makebox(0,0)[lc]{\(2;4\vert 7\)}}
\put(-0.8,-9){\makebox(0,0)[lc]{\(1;5\)}}

\put(1.2,-6.1){\vector(-1,-1){1.5}}
\put(1.3,-6.1){\makebox(0,0)[lc]{last semi-metabelian bifurcation}}

\put(1.2,-8){\makebox(0,0)[lc]{\(2;1\)}}
\put(1.2,-10){\makebox(0,0)[lc]{\(2;8\vert 12\)}}

\put(3.2,-10){\makebox(0,0)[lc]{\(2;4\), periodic root \(T_7\)}}
\put(3.8,-10.3){\makebox(0,0)[lc]{in Formula \eqref{eqn:a1ExcitedD6}}}
\put(3.2,-11){\makebox(0,0)[lc]{\(1;i\)}}

\put(4.2,-11){\makebox(0,0)[lc]{\(1;1\)}}
\put(4.2,-12){\makebox(0,0)[lc]{\(1;i\)}}

\put(5.2,-12){\makebox(0,0)[lc]{\(1;1\)}}
\put(5.2,-13){\makebox(0,0)[lc]{\(1;i\)}}

\put(-9.7,-4){\makebox(0,0)[rc]{\(2;3\)}}
\put(-9.7,-5){\makebox(0,0)[rc]{\(1;1\)}}
\put(-9.7,-7){\makebox(0,0)[rc]{\(2;i\)}}
\put(-9.5,-7.5){\makebox(0,0)[cc]{\(i=4,6\)}}
\put(-9.5,-8){\makebox(0,0)[lc]{\(\mathrm{E}.9\)}}
\put(-9.5,-8.5){\makebox(0,0)[lc]{\((3,3)\)}}

\put(-7.7,-5){\makebox(0,0)[rc]{\(2;1\)}}
\put(-7.7,-6){\makebox(0,0)[rc]{\(1;1\)}}
\put(-7.7,-8){\makebox(0,0)[rc]{\(2;4\)}}
\put(-7.5,-9){\makebox(0,0)[lc]{\(\mathrm{D}.6\)}}
\put(-7.5,-9.5){\makebox(0,0)[lc]{\((9,3)\)}}

\put(-5.7,-6){\makebox(0,0)[rc]{\(2;1\)}}
\put(-5.7,-7){\makebox(0,0)[rc]{\(1;1\)}}
\put(-5.7,-9){\makebox(0,0)[rc]{\(2;4\)}}
\put(-5.5,-10){\makebox(0,0)[lc]{\(\mathrm{D}.6\)}}
\put(-5.5,-10.5){\makebox(0,0)[lc]{\((27,3)\)}}

\put(-3.7,-8){\makebox(0,0)[rc]{\(2;1\)}}
\put(-3.7,-10){\makebox(0,0)[rc]{\(2;4\)}}
\put(-3.5,-11){\makebox(0,0)[lc]{\(\mathrm{D}.6\)}}
\put(-3.5,-11.5){\makebox(0,0)[lc]{\((81,3)\)}}

\put(-2.1,-13.5){\vector(-1,4){1.3}}
\put(-2.3,-13.6){\makebox(0,0)[ct]{last non-metabelian bifurcation}}

\put(-1.7,-10){\makebox(0,0)[rc]{\(2;4\)}}
\put(-1.7,-11){\makebox(0,0)[rc]{\(1;1\)}}
\put(-1.5,-12){\makebox(0,0)[lc]{\(\mathrm{D}.6\)}}
\put(-1.5,-12.5){\makebox(0,0)[lc]{\((243,3)\)}}

\put(1.2,-11){\makebox(0,0)[lc]{\(1;1\)}}
\put(1.2,-12){\makebox(0,0)[lc]{\(1;1\)}}
\put(0.5,-12.5){\makebox(0,0)[lc]{\(i=2,3\)}}
\put(0.5,-13){\makebox(0,0)[lc]{\(\mathrm{D}.6\)}}
\put(0.5,-13.5){\makebox(0,0)[lc]{\((729,3)\)}}

\put(3.2,-12){\makebox(0,0)[lc]{\(1;1\)}}
\put(3.2,-13){\makebox(0,0)[lc]{\(1;1\)}}
\put(2.5,-13.5){\makebox(0,0)[lc]{\(i=2,3\)}}
\put(2.5,-14){\makebox(0,0)[lc]{\(\mathrm{D}.6\)}}
\put(2.5,-14.5){\makebox(0,0)[lc]{\((2187,3)\)}}

\put(4.2,-13){\makebox(0,0)[lc]{\(1;1\)}}
\put(4.2,-14){\makebox(0,0)[lc]{\(1;1\)}}
\put(3.8,-14.5){\makebox(0,0)[lc]{\(i=2,3\)}}
\put(3.5,-15){\makebox(0,0)[lc]{\(\mathrm{D}.6\)}}
\put(3.5,-15.5){\makebox(0,0)[lc]{\((6561,3)\)}}

\put(5.2,-14){\makebox(0,0)[lc]{\(1;1\)}}
\put(5.2,-15){\makebox(0,0)[lc]{\(1;1\)}}
\put(4.8,-15.5){\makebox(0,0)[lc]{\(i=2,3\)}}
\put(4.5,-16){\makebox(0,0)[lc]{\(\mathrm{D}.6\)}}
\put(4.5,-16.5){\makebox(0,0)[lc]{\((19683,3)\)}}

\put(4.8,-9.5){\makebox(0,0)[cc]{main trunk of type \(\mathrm{a}.1\)}}
\put(6.2,-13){\makebox(0,0)[lc]{etc.}}

\put(-4,1.4){\makebox(0,0)[lc]{Legend:}}

\put(-2.5,1.4){\circle{0.2}}
\put(-2.5,1.4){\circle*{0.1}}
\put(-2,1.4){\makebox(0,0)[lc]{\(\ldots\) abelian}}

\put(-2.5,1){\circle{0.2}}
\put(-2,1){\makebox(0,0)[lc]{\(\ldots\) metabelian}}

\put(-2.5,0.6){\circle*{0.2}}
\put(-2,0.6){\makebox(0,0)[lc]{\(\ldots\) metabelian with bifurcation}}

\put(-2.6,0.1){\framebox(0.2,0.2){}}
\put(-2,0.2){\makebox(0,0)[lc]{\(\ldots\) non-metabelian}}

\put(-2.6,-0.3){\framebox(0.2,0.2){}}
\put(-2.5,-0.2){\circle*{0.1}}
\put(-2,-0.2){\makebox(0,0)[lc]{\(\ldots\) non-metabelian with bifurcation}}

\put(-2.6,-0.7){\framebox(0.2,0.2){}}
\put(-2.5,-0.6){\circle{0.1}}
\put(-2,-0.6){\makebox(0,0)[lc]{\(\ldots\) non-metabelian Schur \(\sigma\)}}

\end{picture}

}

\end{figure}


Finally, we come to the first excited state of type \(\mathrm{D}.6\).

\begin{theorem}
\label{thm:a1ChainExcitedD6}
For each logarithmic exponent \(e\ge 7\),
a metabelian \(\mathrm{CF}\)-group \(T_e\)
with commutator quotient \(T_e/T_e^\prime\simeq (3^e,3)\hat{=}(e1)\),
punctured transfer kernel type \(\mathrm{a}.1\), \(\varkappa(T_e)\sim (000;0)\),
and rank distribution \(\varrho\sim (2,2,2;3)\)
is given by the \textbf{periodic sequence} of iterated \(p\)-descendants
\begin{equation}
\label{eqn:a1ExcitedD6}
T_e\simeq\mathrm{SmallGroup}(6561,85)(-\#2;1)^2-\#2;4(-\#1;1)^{e-7}.
\end{equation}
These groups have logarithmic order \(\mathrm{lo}(T_e)=7+e\) and
first abelian quotient invariants \(\alpha_1(T_e)\sim (e1,e1,e1;(e-1)33)\).
They form the infinite main trunk of a descendant tree
with finite double-twigs of depth three.
For each integer \(e\ge 7\), the doublet
\begin{equation}
\label{eqn:ExcitedD6M}
M_{e,i}\simeq\mathrm{SmallGroup}(6561,85)(-\#2;1)^2-\#2;4(-\#1;1)^{e-7}-\#1;i, \qquad i\in\lbrace 2,3\rbrace,
\end{equation}
respectively
\begin{equation}
\label{eqn:ExcitedD6}
G_{e,i}\simeq\mathrm{SmallGroup}(6561,85)(-\#2;1)^2-\#2;4(-\#1;1)^{e-7}-\#1;i(-\#1;1)^2, \qquad i\in\lbrace 2,3\rbrace,
\end{equation}
is the \textbf{unique pair} of \textbf{metabelian} \(\sigma\)-groups \(M_{e,i}\simeq G_{e,i}/G_{e,i}^{\prime\prime}\),
respectively \textbf{non-metabelian} Schur \(\sigma\)-groups \(G_{e,i}\), with
commutator quotient \((3^e,3)\hat{=}(e1)\),
punctured transfer kernel type \(\mathrm{D}.6\), \(\varkappa\sim (123;1)\),
\(\mathrm{lo}(M_{e,i})=8+e\), respectively \(\mathrm{lo}(G_{e,i})=10+e\), and
first abelian quotient invariants \(\alpha_1\sim ((e+1)1,(e+1)1,(e+1)1;e33)\).
See Figure
\ref{fig:SchurSigma7D6}.
\end{theorem}


\begin{proof}
Let \(s_2=\lbrack y,x\rbrack\) denote the main commutator,
\(\forall_{j=3}^7\) \(s_j=\lbrack s_{j-1},x\rbrack\)
and \(t_3=\lbrack s_2,y\rbrack\) higher commutators,
and \(w=x^{3^e}\) the last non-trivial power.
For each \(e\ge 7\), we have parametrized presentations with \(i\in\lbrace 2,3\rbrace\)
\begin{equation}
\label{eqn:Pres7D6}
\begin{aligned}
T_{e+1} &= \langle x,y\mid w^3=1,\ y^3=s_3^2s_4,\ \forall_{j=2}^4\,s_j^3=s_{j+2}^2s_{j+3},\ s_5^3=s_7^2,\ t_3=s_7\rangle, \\
M_{e,i} &= \langle x,y\mid w^3=1,\ y^3=s_3^2s_4,\ \forall_{j=2}^4\,s_j^3=s_{j+2}^2s_{j+3},\ s_5^3=s_7^2,\ t_3=s_7w^{i-1}\rangle.
\end{aligned}
\end{equation}
The \(p\)-class is given by
\(\mathrm{cl}_p(M_{e,i})=e+1\) respectively \(\mathrm{cl}_p(T_{e+1})=e+1\).
For \(e\ge 7\), the last non-trivial lower \(p\)-central is given by
\(P_{e}(M_{e,i})=\langle w\rangle\) respectively \(P_{e}(T_{e+1})=\langle w\rangle\),
whence all three groups share the common \(p\)-parent
\(M_{e,i}/P_{e}(M_{e,i})\simeq T_{e+1}/P_{e}(T_{e+1})\simeq T_e\).

Repeated recursive applications of the \(p\)-group generation algorithm
eventually produce Figure
\ref{fig:SchurSigma7D6},
and thus confirm Formulas
\eqref{eqn:a1ExcitedD6},
\eqref{eqn:ExcitedD6M}
and
\eqref{eqn:ExcitedD6}.
All groups have \(\varrho\sim (2,2,2;3)\).
\end{proof}


\section{Elevated Schur \(\sigma\)-groups \(G\) with \(G/G^\prime\simeq (3^e,3)\), \(e\ge 9\)}
\label{s:NonMetabelian9}

\noindent
Investigation of Schur \(\sigma\)-groups \(G\) with
non-elementary bicyclic commutator quotient \(G/G^\prime\simeq (3^e,3)\), \(e\ge 4\),
punctured transfer kernel type \(\mathrm{B}.18\), \(\varkappa(G)\sim (144;4)\),
and lowest possible logarithmic order \(\mathrm{lo}(G)=19+e\)
is firmly based on a crucial vertex \(T_4\) of the main trunk
with associated scaffold type \(\mathrm{b}.31\), \(\varkappa(T_4)\sim (044;4)\).
The identifier of this fork between non-metabelian Schur \(\sigma\)-groups \(G\)
and their metabelianizations \(M=G/G^{\prime\prime}\) is
\(T_4=\langle 2187,3\rangle-\#3;2\).
It may be called a \textit{bifurcation of infinite order},
since it is responsible for all values \(e\ge 4\) of the logarithmic exponent.
All the groups involved have elevated rank distribution \(\varrho\sim (3,3,3;3)\)
\cite{Ma2021b}.


Periodicity of the metabelianizations \(M\)
sets in for \(e\ge 5\) already.


\begin{theorem}
\label{thm:b31ChainMetabelian}
For each logarithmic exponent \(e\ge 5\),
a metabelian \(\mathrm{BCF}\)-group \(U_e\)
with commutator quotient \(U_e/U_e^\prime\simeq (3^e,3)\hat{=}(e1)\),
punctured transfer kernel type \(\mathrm{b}.31\), \(\varkappa(U_e)\sim (044;4)\),
and elevated rank distribution \(\varrho\sim (3,3,3;3)\)
is given by the \textbf{periodic sequence} of iterated \(p\)-descendants
\begin{equation}
\label{eqn:b31Metabelian}
U_e\simeq\mathrm{SmallGroup}(2187,3)-\#3;2-\#2;93(-\#1;1)^{e-5}.
\end{equation}
These groups have logarithmic order \(\mathrm{lo}(U_e)=7+e\) and
first abelian quotient invariants \(\alpha_1(U_e)\sim (e21,e11,e11;(e-1)21)\).
They form the infinite main trunk of a descendant tree
with finite twigs of depth one,
consisting of a metabelian doublet and \(24\) non-metabelian vertices.
For each integer \(e\ge 5\), the doublet
\begin{equation}
\label{eqn:B18Metabelian}
M_{e,i}\simeq\mathrm{SmallGroup}(2187,3)-\#3;2-\#2;93(-\#1;1)^{e-5}-\#1;i, \qquad i\in\lbrace 2,3\rbrace,
\end{equation}
is the \textbf{unique pair} of \textbf{metabelian} \(\sigma\)-groups with
commutator quotient \((3^e,3)\hat{=}(e1)\),
punctured transfer kernel type \(\mathrm{B}.18\), \(\varkappa\sim (144;4)\),
\(\mathrm{lo}(M_{e,i})=8+e\), and
first abelian quotient invariants \(\alpha_1(M_{e,i})\sim ((e+1)21,e11,e11;(e-1)21)\).
See Figure
\ref{fig:SchurSigma59}.
\end{theorem}


\begin{proof}
Let \(s_2=t_2=\lbrack y,x\rbrack\) denote the main commutator,
\(\forall_{j=3}^5\) \(s_j=\lbrack s_{j-1},x\rbrack\) 
and \(t_j=\lbrack t_{j-1},y\rbrack\) higher commutators,
and \(w=x^{3^e}\) the last non-trivial power.
For each \(e\ge 5\), we have parametrized presentations with \(i\in\lbrace 2,3\rbrace\)
\begin{equation}
\label{eqn:Pres5B18}
\begin{aligned}
U_{e+1} = \langle x,y\mid & w^3=1,\ y^3=s_3s_4^2,\ s_2^3=s_4t_4^2,\ s_3^3=s_5,\ t_3^3=s_5^2, \\
& \lbrack x^3,y\rbrack=s_4t_4s_5^2,\ \lbrack x^3,s_2\rbrack=s_5,\ t_5=s_5\rangle, \\
M_{e,i} = \langle x,y\mid & w^3=1,\ y^3=s_3s_4^2,\ s_2^3=s_4t_4^2w^{i-1},\ s_3^3=s_5,\ t_3^3=s_5^2(w^{i-1})^2, \\
& \lbrack x^3,y\rbrack=s_4t_4s_5^2(w^{i-1})^2,\ \lbrack x^3,s_2\rbrack=s_5,\ t_5=s_5w^{i-1}\rangle.
\end{aligned}
\end{equation}
The \(p\)-class is given by
\(\mathrm{cl}_p(M_{e,i})=e+1\) respectively \(\mathrm{cl}_p(U_{e+1})=e+1\).
For \(e\ge 5\), the last non-trivial lower \(p\)-central is given by
\(P_{e}(M_{e,i})=\langle w\rangle\) respectively \(P_{e}(U_{e+1})=\langle w\rangle\),
whence all three groups share the common \(p\)-parent
\(M_{e,i}/P_{e}(M_{e,i})\simeq U_{e+1}/P_{e}(U_{e+1})\simeq U_e\).

Repeated recursive applications of the \(p\)-group generation algorithm
eventually produce the upper part of Figure
\ref{fig:SchurSigma59},
and thus confirm Formulas
\eqref{eqn:b31Metabelian}
and
\eqref{eqn:B18Metabelian}.
Cf.
\cite[Cor. 1]{Ma2021b}.
\end{proof}


\noindent
Periodicity of the Schur \(\sigma\)-groups \(G\) sets in for \(e\ge 9\).
According to
\cite[Thm. 1--7]{Ma2021b},
the groups with minimal logarithmic order \(\mathrm{lo}(G)=19+e\)
are \(p\)-descendants of the seven roots \(\mathrm{SmallGroup}(2187,3)-\#3;2-\#4;\ell\)
with \(\ell\in\lbrace 24,26,28,30,31,33,37\rbrace\).
We abstain from complete generality
and take \(\mathbf{\ell=37}\), exemplarily.

\begin{theorem}
\label{thm:Elevated37}
A total of \(162\) Schur \(\sigma\)-groups \(G\) with
commutator quotient \(G/G^\prime\simeq (3^e,3)\),
punctured transfer kernel type \(\mathrm{B}.18\), \(\varkappa(G)\sim (144;4)\),
elevated rank distribution  \(\varrho(G)=(3,3,3;3)\),
first abelian quotient invariants \(\alpha_1(G)\sim\lbrack (e+1)21,e11,e11;(e-1)21\rbrack\),
second abelian quotient invariants
\begin{equation}
\label{eqn:Elevated37AQI2}
\begin{aligned}
\alpha_2(G)\sim (e1;&\lbrack (e+1)21;e2111,((e+1)211)^3,((e+1)2)^9 \rbrack, \\
&\lbrack e11;e2111,((e+1)21)^3,((e+1)2)^9 \rbrack, \\
&\lbrack e11;e2111,((e+1)21)^3,((e+1)2)^9 \rbrack; \\
&\lbrack (e-1)21;e2111,(e211)^3,(e21)^8,(e-1)22 \rbrack )
\end{aligned}
\end{equation}
and (minimal) logarithmic order \(\mathrm{lo}(G)=19+e\) is given for each \(e\ge 9\) by the term
\begin{equation}
\label{eqn:Elevated37SchurSigma}
G=T_{9,a,b}\lbrack-\#1;1\rbrack^{e-9}-\#1;i-\#1;k-\#1;1 \text{ with } i\in\lbrace 2,3\rbrace \text{ and } k\in\lbrace 1,2,3\rbrace,
\end{equation}
where \(27\) periodic roots with
\(1\le a\le 9\), \(\tilde{a}=1\) for \(a\in\lbrace 2,6,7\rbrace\), \(\tilde{a}=2\) otherwise, and \(1\le b\le 3\)
are
\begin{equation}
\label{eqn:Elevated37PeriodicRoots}
T_{9,a,b}:=\mathrm{SmallGroup}(2187,3)-\#3;2-\#\mathbf{4;37}-\#3;32-\#4;a-\#2;\tilde{a}-\#2;b.
\end{equation}
\end{theorem}


\begin{proof}
For a fixed step size \(s\ge 1\),
we denote by \(N\) the number of all immediate descendants of a \(3\)-group, and
by \(C\) the number of capable immediate descendants with positive nuclear rank \(\nu\ge 1\).
Generally, let \(X:=\langle 2187,3\rangle-\#3;2-\#4;37-\#3;32\).
This is a non-metabelian \(3\)-group of type \((729,3)\).
We consider a chain of exo-genetic propagations
\cite{Ma2021a}:
\begin{itemize}
\item
\(X\) has \(N=C=27\) for \(s=\nu=4\)
but only the first \(9\) descendants are of type \((2187,3)\).
\item
Each \(X-\#4;a\) with \(1\le a\le 9\)
has \(N=C=6\) for \(s=\nu=2\)
but only the first, resp. second, descendant,
indicated by \(\tilde{a}\in\lbrace 1,2\rbrace\), is of type \((6561,3)\).
\item
Each \(X-\#4;a-\#2,\tilde{a}\) with \(1\le a\le 9\)
has \(N=C=9\) for \(s=\nu=2\)
but only the first \(3\) descendants are of type \((19683,3)\).
\item
Each \(T_{9,a,b}:=X-\#4;a-\#2,\tilde{a}-\#2;b\) with \(1\le a\le 9\) and \(1\le b\le 3\)
has \(6\) Schur \(\sigma\)-descendants
\(T_{9,a,b}\lbrack-\#1;1\rbrack^{e-9}-\#1;i-\#1;k-\#1;1\) with \(i\in\lbrace 2,3\rbrace\) and \(k\in\lbrace 1,2,3\rbrace\),
for each \(e\ge 9\).
\end{itemize}
Together this census yields \(9\cdot 3\cdot 6=162\) Schur \(\sigma\)-groups, for each \(e\ge 9\).
Cf.
\cite[Thm. 6--7]{Ma2021b}.
\end{proof}

\newpage

\noindent
We expand the proof with more details for the
particular instance \(a=1\), \(\tilde{a}=2\) and \(b=1\).

\begin{theorem}
\label{thm:b31Chain}
For each logarithmic exponent \(e\ge 9\),
a metabelian \(\mathrm{CF}\)-group \(T_{e,1,1}\)
with commutator quotient \(T_{e,1,1}/T_{e,1,1}^\prime\simeq (3^e,3)\hat{=}(e1)\),
punctured transfer kernel type \(\mathrm{b}.31\), \(\varkappa(T_{e,1,1})\sim (044;4)\),
and rank distribution \(\varrho(T_{e,1,1})\sim (3,3,3;3)\)
is given by the \textbf{periodic sequence} of iterated \(p\)-descendants
\begin{equation}
\label{eqn:b31}
T_{e,1,1}\simeq T_{9,1,1}(-\#1;1)^{e-9}.
\end{equation}
These groups have logarithmic order \(\mathrm{lo}(T_{e,1,1})=16+e\) and
first abelian quotient invariants \(\alpha_1(T_{e,1,1})\sim (e21,e11,e11;(e-1)21)\).
They form the infinite main trunk of a descendant tree
with finite double-twigs of depth three.
For each integer \(e\ge 9\), the multiplet (doublet or sextet)
\begin{equation}
\label{eqn:B18Alpha2}
X_{e,i}\simeq T_{9,1,1}(-\#1;1)^{e-9}-\#1;i, \qquad i\in\lbrace 2,3\rbrace,
\end{equation}
respectively
\begin{equation}
\label{eqn:B18Alpha}
Y_{e,i,k}\simeq T_{9,1,1}(-\#1;1)^{e-9}-\#1;i-\#1;k, \qquad i\in\lbrace 2,3\rbrace,\ k\in\lbrace 1,2,3\rbrace,
\end{equation}
respectively
\begin{equation}
\label{eqn:B18Schur}
G_{e,i,k}\simeq T_{9,1,1}(-\#1;1)^{e-9}-\#1;i-\#1;k-\#1;1, \qquad i\in\lbrace 2,3\rbrace,\ k\in\lbrace 1,2,3\rbrace,
\end{equation}
is a \textbf{multiplet} of \textbf{non-metabelian} \(\sigma\)-groups
\((X_{e,i})_{i=2,3}\), \((Y_{e,i,k})_{i=2,3;k=1,2,3}\),
respectively \textbf{non-metabelian Schur} \(\sigma\)-groups \((G_{e,i,k})_{i=2,3;k=1,2,3}\), with
commutator quotient \((3^e,3)\hat{=}(e1)\),
punctured transfer kernel type \(\mathrm{B}.18\), \(\varkappa\sim (144;4)\),
\(\mathrm{lo}(X_{e,i})=17+e\), \(\mathrm{lo}(Y_{e,i,k})=18+e\), respectively \(\mathrm{lo}(G_{e,i,k})=19+e\), and
first abelian quotient invariants \(\alpha_1\sim ((e+1)21,e11,e11;(e-1)21)\).
For \(e\ge 9\), the metabelianizations are isomorphic to the groups in Formula
\eqref{eqn:B18Metabelian},
\begin{equation}
\label{eqn:B18Metabelianization}
X_{e,i}/X_{e,i}^{\prime\prime}\simeq Y_{e,i,k}/Y_{e,i,k}^{\prime\prime}\simeq G_{e,i,k}/G_{e,i,k}^{\prime\prime}\simeq M_{e,i},
 \qquad \text{ for all } i\in\lbrace 2,3\rbrace,\ k\in\lbrace 1,2,3\rbrace.
\end{equation}
See Figure
\ref{fig:SchurSigma59}.
The second derived subgroups
are of constant logarithmic order \(\mathrm{lo}=9,10,11\),
in fact, they are abelian of constant type
\begin{equation}
T_{e,1,1}^{\prime\prime}\simeq X_{e,i}^{\prime\prime}\simeq (32211),\qquad Y_{e,i,k}^{\prime\prime}\simeq (322111),\qquad G_{e,i,k}^{\prime\prime}\simeq (332111).
\end{equation}
\end{theorem}


\begin{proof}
Let \(s_2=t_2=\lbrack y,x\rbrack\) denote the main commutator,
\(\forall_{j=3}^9\) \(s_j=\lbrack s_{j-1},x\rbrack\)
and \(\forall_{j=3}^5\) \(t_j=\lbrack t_{j-1},y\rbrack\) higher commutators,
and \(w=x^{3^e}\) the last non-trivial power.
Since the groups are non-metabelian,
additional commutators in the second derived subgroup are required:
\(\forall_{j=5}^7\) \(u_j=\lbrack s_{j-1},y\rbrack\),
\(v_5=\lbrack t_4,x\rbrack\), and \(v_7=\lbrack u_6,x\rbrack\).
The second derived subgroup also contains \(s_6,\ldots,s_9\),
\textbf{but not} \(w\). These facts, together with all relations
(which are only given partially below) immediately lead to the metabelianizations
\begin{equation}
\label{eqn:b31AndB18Metab}
T_{e+1,1,1}/T_{e+1,1,1}^{\prime\prime}\simeq U_{e+1}, \qquad X_{e,i}/X_{e,i}^{\prime\prime}\simeq M_{e,i},
\end{equation}
for \(e\ge 9\) and \(i\in\lbrace 2,3\rbrace\).

For each \(e\ge 9\), we have extensive parametrized presentations with \(i\in\lbrace 2,3\rbrace\),
which we reduce to the decisive relations \textbf{distinct} for \(T_{e+1,1,1}\) and \(X_{e,i}\).
\begin{equation}
\label{eqn:Pres9B18}
\begin{aligned}
T_{e+1,1,1} = \langle x,y\mid & s_2^3=s_4t_4^2s_5u_5v_5^2s_7^2v_7^2s_8^2,\ t_3^3=s_5u_5^2v_5^2s_6v_7^2s_8^2s_9^2, \\
& \lbrack x^3,y\rbrack=s_4t_4s_5v_5^2s_6s_7u_7v_7s_9,\ t_5=s_5^2u_5^2v_5u_6s_8\rangle, \\
X_{e,i} = \langle x,y\mid & s_2^3=s_4t_4^2s_5u_5v_5^2s_7^2v_7^2s_8^2w^{i-1},\ t_3^3=s_5u_5^2v_5^2s_6v_7^2s_8^2s_9^2(w^{i-1})^2, \\
& \lbrack x^3,y\rbrack=s_4t_4s_5v_5^2s_6s_7u_7v_7s_9(w^{i-1})^2,\ t_5=s_5^2u_5^2v_5u_6s_8w^{i-1}\rangle.
\end{aligned}
\end{equation}
The \(p\)-class is given by
\(\mathrm{cl}_p(X_{e,i})=e+1\) respectively \(\mathrm{cl}_p(T_{e+1,1,1})=e+1\).
For \(e\ge 9\), the last non-trivial lower \(p\)-central is given by
\(P_{e}(X_{e,i})=\langle w\rangle\) respectively \(P_{e}(T_{e+1,1,1})=\langle w\rangle\),
whence all three groups share the common \(p\)-parent
\(X_{e,i}/P_{e}(X_{e,i})\simeq T_{e+1,1,1}/P_{e}(T_{e+1,1,1})\simeq T_{e,1,1}\).

Repeated recursive applications of the \(p\)-group generation algorithm
eventually produce the lower part of Figure
\ref{fig:SchurSigma59},
and thus confirm Formulas
\eqref{eqn:b31},
\eqref{eqn:B18Alpha2},
\eqref{eqn:B18Alpha}
and
\eqref{eqn:B18Schur}.
\end{proof}

\newpage

\begin{figure}[ht]
\caption{Schur \(\sigma\)-groups \(G\) with \(\varrho(G)\sim (3,3,3;3)\), \(G/G^\prime\simeq (3^e,3)\), \(2\le e\le 4\)}
\label{fig:SchurSigma}

{\tiny

\setlength{\unitlength}{0.9cm}
\begin{picture}(14,22.5)(-11,-19.5)

\put(-11,2.5){\makebox(0,0)[cb]{order}}

\put(-11,2){\line(0,-1){21}}
\multiput(-11.1,2)(0,-1){22}{\line(1,0){0.2}}


\put(-10.8,2){\makebox(0,0)[lc]{\(3^2\)}}
\put(-10.8,1){\makebox(0,0)[lc]{\(3^3\)}}
\put(-10.8,0){\makebox(0,0)[lc]{\(3^4\)}}
\put(-10.8,-1){\makebox(0,0)[lc]{\(3^5\)}}
\put(-10.8,-2){\makebox(0,0)[lc]{\(3^6\)}}
\put(-10.8,-3){\makebox(0,0)[lc]{\(3^7\)}}
\put(-10.8,-4){\makebox(0,0)[lc]{\(3^8\)}}
\put(-10.8,-5){\makebox(0,0)[lc]{\(3^9\)}}
\put(-10.8,-6){\makebox(0,0)[lc]{\(3^{10}\)}}
\put(-10.8,-7){\makebox(0,0)[lc]{\(3^{11}\)}}
\put(-10.8,-8){\makebox(0,0)[lc]{\(3^{12}\)}}
\put(-10.8,-9){\makebox(0,0)[lc]{\(3^{13}\)}}
\put(-10.8,-10){\makebox(0,0)[lc]{\(3^{14}\)}}
\put(-10.8,-11){\makebox(0,0)[lc]{\(3^{15}\)}}
\put(-10.8,-12){\makebox(0,0)[lc]{\(3^{16}\)}}
\put(-10.8,-13){\makebox(0,0)[lc]{\(3^{17}\)}}
\put(-10.8,-14){\makebox(0,0)[lc]{\(3^{18}\)}}
\put(-10.8,-15){\makebox(0,0)[lc]{\(3^{19}\)}}
\put(-10.8,-16){\makebox(0,0)[lc]{\(3^{20}\)}}
\put(-10.8,-17){\makebox(0,0)[lc]{\(3^{21}\)}}
\put(-10.8,-18){\makebox(0,0)[lc]{\(3^{22}\)}}
\put(-10.8,-19){\makebox(0,0)[lc]{\(3^{23}\)}}


\put(-11,-19){\vector(0,-1){1}}

\put(-9.1,1.9){\framebox(0.2,0.2){}}
\put(-9,2){\circle*{0.1}}

\put(-9,1){\circle{0.2}}
\put(-9,-1){\circle{0.2}}
\put(-9,-3){\circle*{0.2}}
\put(-8.5,-5){\circle{0.2}}

\put(-7,0){\circle{0.2}}
\put(-7,-2){\circle{0.2}}
\put(-7,-4){\circle*{0.2}}
\put(-6.5,-6){\circle{0.2}}

\put(-5,-3){\circle{0.2}}
\put(-5,-5){\circle*{0.2}}
\put(-4.5,-7){\circle{0.2}}

\put(-3,-6){\circle*{0.2}}
\put(-2.5,-8){\circle{0.2}}

\put(-9.1,-7.1){\framebox(0.2,0.2){}}
\put(-9.1,-9.1){\framebox(0.2,0.2){}}
\put(-9.1,-13.1){\framebox(0.2,0.2){}}
\put(-9.1,-14.1){\framebox(0.2,0.2){}}
\put(-9.1,-16.1){\framebox(0.2,0.2){}}
\put(-9,-16){\circle{0.1}}

\put(-7.1,-8.1){\framebox(0.2,0.2){}}
\put(-7.1,-10.1){\framebox(0.2,0.2){}}
\put(-7.1,-14.1){\framebox(0.2,0.2){}}
\put(-7.1,-15.1){\framebox(0.2,0.2){}}
\put(-7.1,-17.1){\framebox(0.2,0.2){}}
\put(-7,-17){\circle{0.1}}

\put(-5.1,-9.1){\framebox(0.2,0.2){}}
\put(-5.1,-11.1){\framebox(0.2,0.2){}}
\put(-5.1,-15.1){\framebox(0.2,0.2){}}
\put(-5.1,-16.1){\framebox(0.2,0.2){}}
\put(-5.1,-18.1){\framebox(0.2,0.2){}}
\put(-5,-18){\circle{0.1}}

\put(-3.1,-10.1){\framebox(0.2,0.2){}}
\put(-3.1,-12.1){\framebox(0.2,0.2){}}
\put(-3.1,-16.1){\framebox(0.2,0.2){}}
\put(-3.1,-17.1){\framebox(0.2,0.2){}}
\put(-3.1,-19.1){\framebox(0.2,0.2){}}
\put(-3,-19){\circle{0.1}}



\put(-9,2){\line(0,-1){1}}
\put(-9,1){\line(0,-1){2}}
\put(-9,-1){\line(0,-1){2}}
\put(-9,-3){\line(0,-1){4}}
\put(-9,-3){\line(1,-4){0.5}}
\put(-9,-7){\line(0,-1){2}}
\put(-9,-9){\line(0,-1){4}}
\put(-9,-13){\line(0,-1){1}}
\put(-9,-14){\line(0,-1){2}}

\put(-9,2){\line(1,-1){2}}
\put(-7,0){\line(0,-1){2}}
\put(-7,-2){\line(0,-1){2}}
\put(-7,-4){\line(0,-1){4}}
\put(-7,-4){\line(1,-4){0.5}}
\put(-7,-8){\line(0,-1){2}}
\put(-7,-10){\line(0,-1){4}}
\put(-7,-14){\line(0,-1){1}}
\put(-7,-15){\line(0,-1){2}}

\put(-7,0){\line(2,-3){2}}
\put(-5,-3){\line(0,-1){2}}
\put(-5,-5){\line(0,-1){4}}
\put(-5,-5){\line(1,-4){0.5}}
\put(-5,-9){\line(0,-1){2}}
\put(-5,-11){\line(0,-1){4}}
\put(-5,-15){\line(0,-1){1}}
\put(-5,-16){\line(0,-1){2}}

\put(-5,-3){\line(2,-3){2}}
\put(-3,-6){\line(0,-1){4}}
\put(-3,-6){\line(1,-4){0.5}}
\put(-3,-10){\line(0,-1){2}}
\put(-3,-12){\line(0,-1){4}}
\put(-3,-16){\line(0,-1){1}}
\put(-3,-17){\line(0,-1){2}}

\put(-3,-6){\line(1,-1){2}}
\put(-0.8,-8){\makebox(0,0)[lc]{\((243,3)\), continued in Fig. \ref{fig:SchurSigma59}}}


\put(-8.8,2){\makebox(0,0)[lc]{\(\langle 2\rangle\)}}

\put(-8.8,1){\makebox(0,0)[lc]{\(\langle 3\rangle\)}}
\put(-8.8,-1){\makebox(0,0)[lc]{\(\langle 3\rangle\)}}
\put(-8.8,-3){\makebox(0,0)[lc]{\(\langle 64\rangle\)}}
\put(-8.3,-5){\makebox(0,0)[lc]{\(2;36\)}}

\put(-6.8,0){\makebox(0,0)[lc]{\(\langle 3\rangle\)}}
\put(-6.8,-2){\makebox(0,0)[lc]{\(\langle 10\rangle\)}}
\put(-6.8,-4){\makebox(0,0)[lc]{\(\langle 165\rangle\)}}
\put(-6.3,-6){\makebox(0,0)[lc]{\(2;85\)}}

\put(-4.8,-3){\makebox(0,0)[lc]{\(\langle 3\rangle\)}}
\put(-4.8,-5){\makebox(0,0)[lc]{\(2;10\)}}
\put(-4.3,-7){\makebox(0,0)[lc]{\(2;88\)}}

\put(-2.8,-6){\makebox(0,0)[lc]{\(3;2\), the last bifurcation, \(\mathrm{cl}=4=e\)}}
\put(-2.3,-8){\makebox(0,0)[lc]{\(2;100\)}}

\put(-8.8,-7){\makebox(0,0)[lc]{\(4;144\)}}
\put(-8.8,-9){\makebox(0,0)[lc]{\(2;41\)}}
\put(-8.8,-13){\makebox(0,0)[lc]{\(4;1\)}}
\put(-8.8,-14){\makebox(0,0)[lc]{\(1;1\)}}
\put(-8.8,-16){\makebox(0,0)[lc]{\(2;1\)}}
\put(-9,-16.5){\makebox(0,0)[lc]{\(\mathrm{F}.11\)}}
\put(-9,-17){\makebox(0,0)[lc]{\((3,3)\)}}

\put(-6.8,-8){\makebox(0,0)[lc]{\(4;29\)}}
\put(-6.8,-10){\makebox(0,0)[lc]{\(2;41\)}}
\put(-6.8,-14){\makebox(0,0)[lc]{\(4;1\)}}
\put(-6.8,-15){\makebox(0,0)[lc]{\(1;2\)}}
\put(-6.8,-17){\makebox(0,0)[lc]{\(2;1\)}}
\put(-7,-17.5){\makebox(0,0)[lc]{\(\mathrm{B}.18\)}}
\put(-7,-18){\makebox(0,0)[lc]{\((9,3)\)}}

\put(-4.8,-9){\makebox(0,0)[lc]{\(4;51\)}}
\put(-4.8,-11){\makebox(0,0)[lc]{\(2;32\)}}
\put(-4.8,-15){\makebox(0,0)[lc]{\(4;1\)}}
\put(-4.8,-16){\makebox(0,0)[lc]{\(1;1\)}}
\put(-4.8,-18){\makebox(0,0)[lc]{\(2;1\)}}
\put(-5,-18.5){\makebox(0,0)[lc]{\(\mathrm{B}.18\)}}
\put(-5,-19){\makebox(0,0)[lc]{\((27,3)\)}}

\put(-2.8,-10){\makebox(0,0)[lc]{\(4;81\)}}
\put(-2.8,-12){\makebox(0,0)[lc]{\(2;5\)}}
\put(-2.8,-16){\makebox(0,0)[lc]{\(4;1\)}}
\put(-2.8,-17){\makebox(0,0)[lc]{\(1;1\)}}
\put(-2.8,-19){\makebox(0,0)[lc]{\(2;1\)}}
\put(-3,-19.5){\makebox(0,0)[lc]{\(\mathrm{B}.18\)}}
\put(-3,-20){\makebox(0,0)[lc]{\((81,3)\)}}


\put(-4,-0.6){\makebox(0,0)[lc]{Legend:}}

\put(-2.6,-0.7){\framebox(0.2,0.2){}}
\put(-2.5,-0.6){\circle*{0.1}}
\put(-2,-0.6){\makebox(0,0)[lc]{\(\ldots\) abelian}}

\put(-2.5,-1){\circle{0.2}}
\put(-2,-1){\makebox(0,0)[lc]{\(\ldots\) metabelian}}

\put(-2.5,-1.4){\circle*{0.2}}
\put(-2,-1.4){\makebox(0,0)[lc]{\(\ldots\) metabelian with bifurcation}}

\put(-2.6,-1.9){\framebox(0.2,0.2){}}
\put(-2,-1.8){\makebox(0,0)[lc]{\(\ldots\) non-metabelian}}

\put(-2.6,-2.3){\framebox(0.2,0.2){}}
\put(-2.5,-2.2){\circle{0.1}}
\put(-2,-2.2){\makebox(0,0)[lc]{\(\ldots\) non-metabelian Schur \(\sigma\)}}

\end{picture}

}

\end{figure}
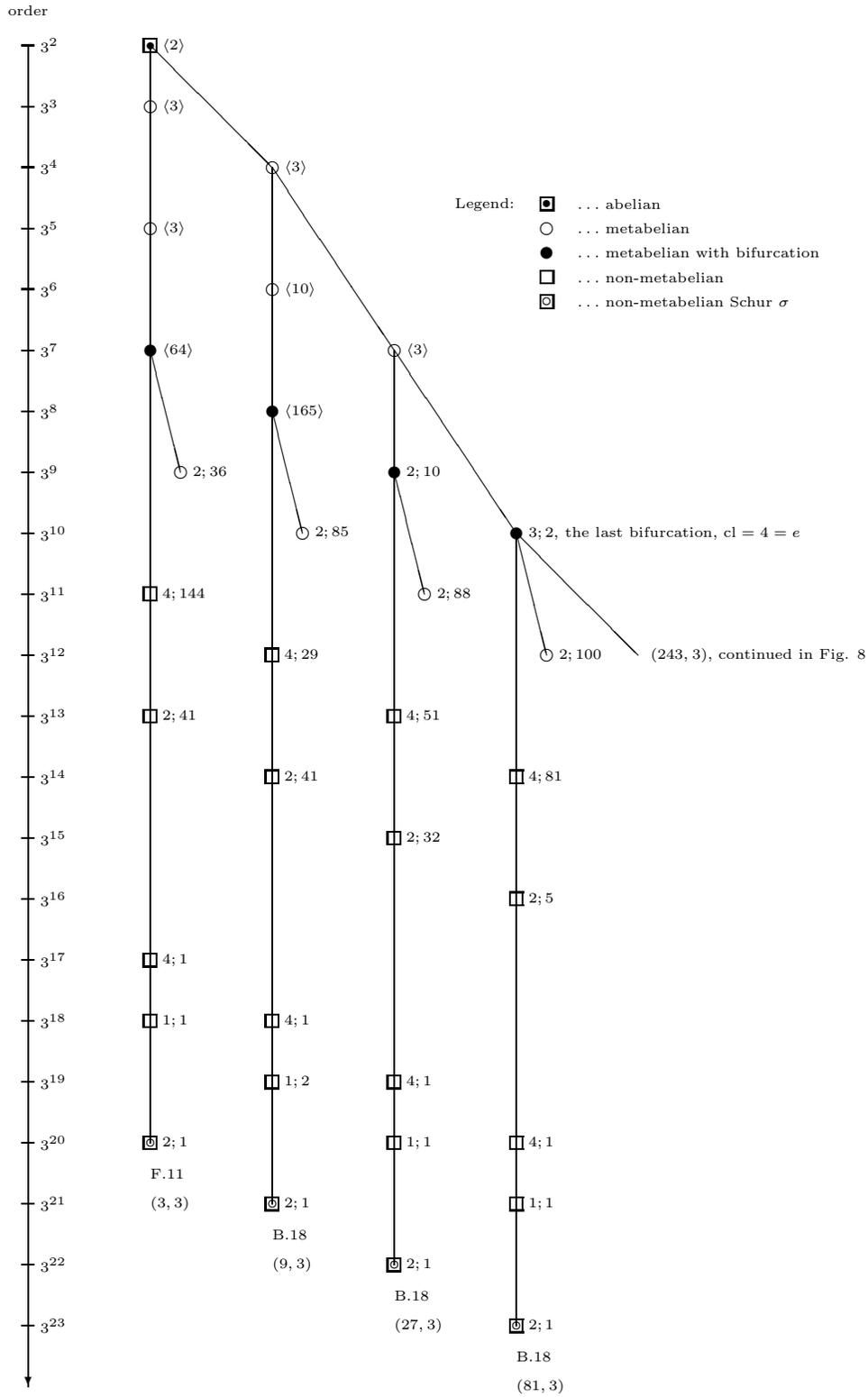

\newpage

\begin{figure}[ht]
\caption{Schur \(\sigma\)-groups \(G\) with \(\varrho(G)\sim (3,3,3;3)\), \(G/G^\prime\simeq (3^e,3)\), \(4\le e\le 13\)}
\label{fig:SchurSigma59}

{\tiny

\setlength{\unitlength}{0.9cm}
\begin{picture}(14,21.5)(-10,-18.5)

\put(-11,2.5){\makebox(0,0)[cb]{order}}

\put(-11,2){\line(0,-1){22}}
\multiput(-11.1,2)(0,-1){23}{\line(1,0){0.2}}


\put(-10.8,2){\makebox(0,0)[lc]{\(3^{10}\)}}
\put(-10.8,1){\makebox(0,0)[lc]{\(3^{11}\)}}
\put(-10.8,0){\makebox(0,0)[lc]{\(3^{12}\)}}
\put(-10.8,-1){\makebox(0,0)[lc]{\(3^{13}\)}}
\put(-10.8,-2){\makebox(0,0)[lc]{\(3^{14}\)}}
\put(-10.8,-3){\makebox(0,0)[lc]{\(3^{15}\)}}
\put(-10.8,-4){\makebox(0,0)[lc]{\(3^{16}\)}}
\put(-10.8,-5){\makebox(0,0)[lc]{\(3^{17}\)}}
\put(-10.8,-6){\makebox(0,0)[lc]{\(3^{18}\)}}
\put(-10.8,-7){\makebox(0,0)[lc]{\(3^{19}\)}}
\put(-10.8,-8){\makebox(0,0)[lc]{\(3^{20}\)}}
\put(-10.8,-9){\makebox(0,0)[lc]{\(3^{21}\)}}
\put(-10.8,-10){\makebox(0,0)[lc]{\(3^{22}\)}}
\put(-10.8,-11){\makebox(0,0)[lc]{\(3^{23}\)}}
\put(-10.8,-12){\makebox(0,0)[lc]{\(3^{24}\)}}
\put(-10.8,-13){\makebox(0,0)[lc]{\(3^{25}\)}}
\put(-10.8,-14){\makebox(0,0)[lc]{\(3^{26}\)}}
\put(-10.8,-15){\makebox(0,0)[lc]{\(3^{27}\)}}
\put(-10.8,-16){\makebox(0,0)[lc]{\(3^{28}\)}}
\put(-10.8,-17){\makebox(0,0)[lc]{\(3^{29}\)}}
\put(-10.8,-18){\makebox(0,0)[lc]{\(3^{30}\)}}
\put(-10.8,-19){\makebox(0,0)[lc]{\(3^{31}\)}}
\put(-10.8,-20){\makebox(0,0)[lc]{\(3^{32}\)}}


\put(-11,-20){\vector(0,-1){1}}


\put(-10,2){\circle*{0.2}}
\put(-9.5,0){\circle{0.2}}
\put(-8,0){\circle{0.2}}
\put(-8,-1){\circle{0.2}}
\put(-7,-1){\circle{0.2}}
\put(-7,-2){\circle{0.2}}


\put(-10.1,-2.1){\framebox(0.2,0.2){}}
\put(-10.1,-4.1){\framebox(0.2,0.2){}}
\put(-10.1,-8.1){\framebox(0.2,0.2){}}
\put(-10.1,-9.1){\framebox(0.2,0.2){}}
\put(-10.1,-11.1){\framebox(0.2,0.2){}}
\put(-10,-11){\circle{0.1}}

\put(-5.1,-2.1){\framebox(0.2,0.2){}}
\put(-5.1,-5.1){\framebox(0.2,0.2){}}
\put(-5.1,-9.1){\framebox(0.2,0.2){}}
\put(-5.1,-10.1){\framebox(0.2,0.2){}}
\put(-5.1,-12.1){\framebox(0.2,0.2){}}
\put(-5,-12){\circle{0.1}}

\put(-4.1,-5.1){\framebox(0.2,0.2){}}
\put(-4.1,-9.1){\framebox(0.2,0.2){}}
\put(-4.1,-11.1){\framebox(0.2,0.2){}}
\put(-4.1,-13.1){\framebox(0.2,0.2){}}
\put(-4,-13){\circle{0.1}}

\put(-3.1,-9.1){\framebox(0.2,0.2){}}
\put(-3.1,-11.1){\framebox(0.2,0.2){}}
\put(-3.1,-13.1){\framebox(0.2,0.2){}}
\put(-3.1,-14.1){\framebox(0.2,0.2){}}
\put(-3,-14){\circle{0.1}}

\put(-2.1,-11.1){\framebox(0.2,0.2){}}
\put(-2.1,-13.1){\framebox(0.2,0.2){}}
\put(-2.1,-14.1){\framebox(0.2,0.2){}}
\put(-2.1,-15.1){\framebox(0.2,0.2){}}
\put(-2,-15){\circle{0.1}}

\put(-1.1,-13.1){\framebox(0.2,0.2){}}
\put(-1.1,-14.1){\framebox(0.2,0.2){}}
\put(-1.1,-15.1){\framebox(0.2,0.2){}}
\put(-1.1,-16.1){\framebox(0.2,0.2){}}
\put(-1,-16){\circle{0.1}}

\put(-0.1,-14.1){\framebox(0.2,0.2){}}
\put(-0.1,-15.1){\framebox(0.2,0.2){}}
\put(-0.1,-16.1){\framebox(0.2,0.2){}}
\put(-0.1,-17.1){\framebox(0.2,0.2){}}
\put(0,-17){\circle{0.1}}

\put(0.9,-15.1){\framebox(0.2,0.2){}}
\put(0.9,-16.1){\framebox(0.2,0.2){}}
\put(0.9,-17.1){\framebox(0.2,0.2){}}
\put(0.9,-18.1){\framebox(0.2,0.2){}}
\put(1,-18){\circle{0.1}}

\put(1.9,-16.1){\framebox(0.2,0.2){}}
\put(1.9,-17.1){\framebox(0.2,0.2){}}
\put(1.9,-18.1){\framebox(0.2,0.2){}}
\put(1.9,-19.1){\framebox(0.2,0.2){}}
\put(2,-19){\circle{0.1}}

\put(2.9,-17.1){\framebox(0.2,0.2){}}
\put(2.9,-18.1){\framebox(0.2,0.2){}}
\put(2.9,-19.1){\framebox(0.2,0.2){}}
\put(2.9,-20.1){\framebox(0.2,0.2){}}
\put(3,-20){\circle{0.1}}


\put(-10,2){\line(0,-1){4}}
\put(-10,2){\line(1,-4){0.5}}
\put(-10,-2){\line(0,-1){2}}
\put(-10,-4){\line(0,-1){4}}
\put(-10,-8){\line(0,-1){1}}
\put(-10,-9){\line(0,-1){2}}

\put(-10,2){\line(1,-1){2}}
\put(-8,0){\line(0,-1){1}}
\put(-8,0){\line(1,-1){1}}
\put(-7,-1){\line(0,-1){1}}
\put(-7,-1){\line(1,-1){0.7}}

\put(-10,2){\line(5,-4){5}}
\put(-5,-2){\line(0,-1){3}}
\put(-5,-5){\line(0,-1){4}}
\put(-5,-9){\line(0,-1){1}}
\put(-5,-10){\line(0,-1){2}}

\put(-5,-2){\line(1,-3){1}}
\put(-4,-5){\line(0,-1){4}}
\put(-4,-9){\line(0,-1){2}}
\put(-4,-11){\line(0,-1){2}}

\put(-4,-5){\line(1,-4){1}}
\put(-3,-9){\line(0,-1){2}}
\put(-3,-11){\line(0,-1){2}}
\put(-3,-13){\line(0,-1){1}}

\put(-3,-9){\line(1,-2){1}}
\put(-2,-11){\line(0,-1){2}}
\put(-2,-13){\line(0,-1){1}}
\put(-2,-14){\line(0,-1){1}}

\put(-2,-11){\line(1,-2){1}}
\put(-1,-13){\line(0,-1){1}}
\put(-1,-14){\line(0,-1){1}}
\put(-1,-15){\line(0,-1){1}}

\put(-1,-13){\line(1,-1){1}}
\put(0,-14){\line(0,-1){1}}
\put(0,-15){\line(0,-1){1}}
\put(0,-16){\line(0,-1){1}}

\put(0,-14){\line(1,-1){1}}
\put(1,-15){\line(0,-1){1}}
\put(1,-16){\line(0,-1){1}}
\put(1,-17){\line(0,-1){1}}

\put(1,-15){\line(1,-1){1}}
\put(2,-16){\line(0,-1){1}}
\put(2,-17){\line(0,-1){1}}
\put(2,-18){\line(0,-1){1}}

\put(2,-16){\line(1,-1){1}}
\put(3,-17){\line(0,-1){1}}
\put(3,-18){\line(0,-1){1}}
\put(3,-19){\line(0,-1){1}}

\put(3,-17){\line(1,-1){1}}



\put(-9.8,2){\makebox(0,0)[lc]{\(3;2\), the last bifurcation, \(\mathrm{cl}=4=e\), continued from Fig. \ref{fig:SchurSigma}}}
\put(-9.3,0){\makebox(0,0)[lc]{\(2;100\)}}
\put(-9.5,-0.3){\makebox(0,0)[lc]{\((81,3)\)}}
\put(-7.8,0){\makebox(0,0)[lc]{\(2;93\), metabelian periodic root \(U_5\)}}
\put(-7,-0.3){\makebox(0,0)[lc]{in Formula \eqref{eqn:b31Metabelian}}}

\put(-7.8,-1){\makebox(0,0)[lc]{\(1;i\)}}
\put(-8,-1.3){\makebox(0,0)[cc]{\(i=2,3\)}}
\put(-8,-1.6){\makebox(0,0)[cc]{\(\mathrm{B}.18\)}}
\put(-8,-1.9){\makebox(0,0)[cc]{\((243,3)\)}}

\put(-6.8,-1){\makebox(0,0)[lc]{\(1;1\)}}

\put(-6.8,-2){\makebox(0,0)[lc]{\(1;i\)}}
\put(-7,-2.3){\makebox(0,0)[cc]{\(i=2,3\)}}
\put(-7,-2.6){\makebox(0,0)[cc]{\(\mathrm{B}.18\)}}
\put(-7,-2.9){\makebox(0,0)[cc]{\((729,3)\)}}

\put(-6.2,-1.7){\makebox(0,0)[lc]{etc.}}


\put(-9.8,-2){\makebox(0,0)[lc]{\(4;81\)}}
\put(-9.8,-4){\makebox(0,0)[lc]{\(2;5\)}}
\put(-9.8,-8){\makebox(0,0)[lc]{\(4;1\)}}
\put(-9.8,-9){\makebox(0,0)[lc]{\(1;1\)}}
\put(-9.8,-11){\makebox(0,0)[lc]{\(2;1\)}}
\put(-10,-11.5){\makebox(0,0)[lc]{\(\mathrm{B}.18\)}}
\put(-10,-12){\makebox(0,0)[lc]{\((81,3)\)}}

\put(-4.8,-2){\makebox(0,0)[lc]{\(4;37\), exemplary selection}}
\put(-4.8,-5){\makebox(0,0)[lc]{\(3;73\)}}
\put(-4.8,-9){\makebox(0,0)[lc]{\(4;1\)}}
\put(-4.8,-10){\makebox(0,0)[lc]{\(1;2\)}}
\put(-4.8,-12){\makebox(0,0)[lc]{\(2;1\)}}
\put(-5,-12.5){\makebox(0,0)[rc]{\(\mathrm{B}.18\)}}
\put(-5,-13){\makebox(0,0)[rc]{\((243,3)\)}}

\put(-3.8,-5){\makebox(0,0)[lc]{\(3;32\)}}
\put(-3.8,-9){\makebox(0,0)[lc]{\(4;10\)}}
\put(-3.8,-11){\makebox(0,0)[lc]{\(2;2\)}}
\put(-3.8,-13){\makebox(0,0)[lc]{\(2;1\)}}
\put(-4,-13.5){\makebox(0,0)[rc]{\(\mathrm{B}.18\)}}
\put(-4,-14){\makebox(0,0)[rc]{\((729,3)\)}}

\put(-2.8,-9){\makebox(0,0)[lc]{\(4;a\) with \(a=1\)}}
\put(-2.8,-11){\makebox(0,0)[lc]{\(2;4\)}}
\put(-2.8,-13){\makebox(0,0)[lc]{\(2;1\)}}
\put(-2.8,-14){\makebox(0,0)[lc]{\(1;1\)}}
\put(-3,-14.5){\makebox(0,0)[rc]{\(\mathrm{B}.18\)}}
\put(-3,-15){\makebox(0,0)[rc]{\((2187,3)\)}}

\put(-1.8,-11){\makebox(0,0)[lc]{\(2;\tilde{a}\) with \(\tilde{a}=2\)}}
\put(-1.8,-13){\makebox(0,0)[lc]{\(2;4\)}}
\put(-1.8,-14){\makebox(0,0)[lc]{\(1;1\)}}
\put(-1.8,-15){\makebox(0,0)[lc]{\(1;1\)}}
\put(-2,-15.5){\makebox(0,0)[rc]{\(\mathrm{B}.18\)}}
\put(-2,-16){\makebox(0,0)[rc]{\((6561,3)\)}}

\put(-0.8,-13){\makebox(0,0)[lc]{\(2;b\) with \(b=1\), non-metabelian periodic root \(T_{9,1,1}\)}}
\put(1.2,-13.3){\makebox(0,0)[lc]{in Formula \eqref{eqn:Elevated37PeriodicRoots} and \eqref{eqn:b31}}}
\put(-0.8,-14){\makebox(0,0)[lc]{\(1;i\)}}
\put(-0.8,-15){\makebox(0,0)[lc]{\(1;k\)}}
\put(-0.8,-16){\makebox(0,0)[lc]{\(1;1\)}}
\put(-1,-16.4){\makebox(0,0)[rc]{\(i=2,3\)}}
\put(-1,-16.7){\makebox(0,0)[rc]{\(\mathrm{B}.18\)}}
\put(-1,-17){\makebox(0,0)[rc]{\((19683,3)\)}}

\put(0.2,-14){\makebox(0,0)[lc]{\(1;1\)}}
\put(0.2,-15){\makebox(0,0)[lc]{\(1;i\)}}
\put(0.2,-16){\makebox(0,0)[lc]{\(1;k\)}}
\put(0.2,-17){\makebox(0,0)[lc]{\(1;1\)}}
\put(0,-17.4){\makebox(0,0)[rc]{\(i=2,3\)}}
\put(0,-17.7){\makebox(0,0)[rc]{\(\mathrm{B}.18\)}}
\put(0,-18){\makebox(0,0)[rc]{\((59049,3)\)}}

\put(1.2,-15){\makebox(0,0)[lc]{\(1;1\)}}
\put(1.2,-16){\makebox(0,0)[lc]{\(1;i\)}}
\put(1.2,-17){\makebox(0,0)[lc]{\(1;k\)}}
\put(1.2,-18){\makebox(0,0)[lc]{\(1;1\)}}
\put(1,-18.4){\makebox(0,0)[rc]{\(i=2,3\)}}
\put(1,-18.7){\makebox(0,0)[rc]{\(\mathrm{B}.18\)}}
\put(1,-19){\makebox(0,0)[rc]{\((177147,3)\)}}

\put(2.2,-16){\makebox(0,0)[lc]{\(1;1\)}}
\put(2.2,-17){\makebox(0,0)[lc]{\(1;i\)}}
\put(2.2,-18){\makebox(0,0)[lc]{\(1;k\)}}
\put(2.2,-19){\makebox(0,0)[lc]{\(1;1\)}}
\put(2,-19.4){\makebox(0,0)[rc]{\(i=2,3\)}}
\put(2,-19.7){\makebox(0,0)[rc]{\(\mathrm{B}.18\)}}
\put(2,-20){\makebox(0,0)[rc]{\((531441,3)\)}}

\put(3.2,-17){\makebox(0,0)[lc]{\(1;1\)}}
\put(3.2,-18){\makebox(0,0)[lc]{\(1;i\)}}
\put(3.2,-19){\makebox(0,0)[lc]{\(1;k\)}}
\put(3.2,-20){\makebox(0,0)[lc]{\(1;1\)}}
\put(3,-20.4){\makebox(0,0)[rc]{\(i=2,3\)}}
\put(3,-20.7){\makebox(0,0)[rc]{\(\mathrm{B}.18\)}}
\put(3,-21){\makebox(0,0)[rc]{\((1594323,3)\)}}

\put(0,-20.4){\makebox(0,0)[rc]{\(k=1,2,3\)}}
\put(1.7,-15.5){\makebox(0,0)[lc]{main trunk of type \(\mathrm{b}.31\)}}
\put(4.2,-18){\makebox(0,0)[lc]{etc.}}


\put(-4,1.4){\makebox(0,0)[lc]{Legend:}}

\put(-2.5,1.4){\circle{0.2}}
\put(-2,1.4){\makebox(0,0)[lc]{\(\ldots\) metabelian}}

\put(-2.5,1){\circle*{0.2}}
\put(-2,1){\makebox(0,0)[lc]{\(\ldots\) metabelian with bifurcation}}

\put(-2.6,0.5){\framebox(0.2,0.2){}}
\put(-2,0.6){\makebox(0,0)[lc]{\(\ldots\) non-metabelian}}

\put(-2.6,0.1){\framebox(0.2,0.2){}}
\put(-2.5,0.2){\circle{0.1}}
\put(-2,0.2){\makebox(0,0)[lc]{\(\ldots\) non-metabelian Schur \(\sigma\)}}

\end{picture}

}

\end{figure}

\newpage

\section{Acknowledgements}
\label{s:Gratifications}

\noindent
The author acknowledges
support by the Austrian Science Fund (FWF): Projects J0497-PHY and P26008-N25,
and by the Research Executive Agency of the European Union (EUREA).




\begin{thebibliography}{XX}
%
\bibitem{Ag1998}
M. Arrigoni,
\textit{On Schur \(\sigma\)-groups},
Math. Nachr.
\textbf{192}
(1998),
71--89.
%
\bibitem{Ar1927}
E. Artin,
\textit{Beweis des allgemeinen Reziprozit\"atsgesetzes},
Abh. Math. Sem. Univ. Hamburg
\textbf{5}
(1927),
353--363.
%
\bibitem{Ar1929}
E. Artin,
\textit{Idealklassen in Oberk\"orpern und allgemeines Reziprozit\"atsgesetz},
Abh. Math. Sem. Univ. Hamburg
\textbf{7}
(1929),
46--51.
%
\bibitem{AHL1977}
J. A. Ascione, G. Havas, and C. R. Leedham-Green,
\textit{A computer aided classification of certain groups of prime power order},
Bull. Austral. Math. Soc.
\textbf{17}
(1977),
257--274.
%
%
\bibitem{BEO2005}
H. U. Besche, B. Eick, and E. A. O'Brien,
\textit{The SmallGroups Library --- a Library of Groups of Small Order},
2005,
an accepted and refereed GAP package, available also in MAGMA.
%
\bibitem{BCP1997}
W. Bosma, J. Cannon, and C. Playoust,
\textit{The Magma algebra system. I. The user language}, 
J. Symbolic Comput.
\textbf{24}
(1997),
235--265.
%
\bibitem{BCFS2021}
W. Bosma, J. J. Cannon, C. Fieker, A. Steels (eds.),
\textit{Handbook of Magma functions},
Ed. 2.26,
Sydney,
2021.
%
\bibitem{BBH2017}
N. Boston, M. R. Bush and F. Hajir,
\textit{Heuristics for \(p\)-class towers of imaginary quadratic fields},
Math. Ann.
\textbf{368}
(2017),
No. 1,
633--669,
DOI 10.1007/s00208-016-1449-3.
%
\bibitem{BuMa2015}
M. R. Bush and D. C. Mayer,
\textit{\(3\)-class field towers of exact length \(3\)},
J. Number Theory
\textbf{147}
(2015),
766--777, \\
DOI 10.1016/j.jnt.2014.08.010.
%
\bibitem{ELNO2013}
B. Eick, C. R. Leedham-Green, M. F. Newman, and E. A. O'Brien,
\textit{On the classification of groups of prime-power order by coclass:
The \(3\)-groups of coclass \(2\)},
Int. J. Algebra Comput.
\textbf{23}
(2013),
1243--1288.
%
\bibitem{Fi2001}
C. Fieker,
\textit{Computing class fields via the Artin map},
Math. Comp.
\textbf{70}
(2001),
No. 235,
1293--1303.
%
\bibitem{GNO2006}
G. Gamble, W. Nickel, and E. A. O'Brien,
\textit{ANU \(p\)-Quotient --- \(p\)-Quotient and \(p\)-Group Generation Algorithms},
2006,
an accepted GAP package, available also in MAGMA.
%
%
\bibitem{HEO2005}
D. F. Holt, B. Eick, and E. A. O'Brien,
\textit{Handbook of computational group theory},
Discrete mathematics and its applications,
Chapman and Hall/CRC Press,
Boca Raton,
2005.
%
\bibitem{KoVe1975}
H. Koch und B. B. Venkov,
\textit{\"Uber den \(p\)-Klassenk\"orperturm eines imagin\"ar-quadra\-tischen Zahlk\"orpers},
Ast\'erisque
\textbf{24--25}
(1975),
57--67.
%
\bibitem{MAGMA2021}
MAGMA Developer Group,
MAGMA \textit{Computational Algebra System},
Version 2.26-8,
Univ. Sydney,
2021, \\
\texttt{(http://magma.maths.usyd.edu.au)}.
%
%
\bibitem{Ma2012a}
D. C. Mayer,
\textit{Transfers of metabelian \(p\)-groups},
Monatsh. Math.
\textbf{166}
(2012),
No. 3--4,
467--495, \\
DOI 10.1007/s00605-010-0277-x.
%
\bibitem{Ma2015a}
D. C. Mayer,
\textit{Periodic bifurcations in descendant trees of finite \(p\)-groups},
Adv. Pure Math.
\textbf{5}
(2015),
No. 1,
162--195,
DOI 10.4236/apm.2015.54020.
%
\bibitem{Ma2015b}
D. C. Mayer,
\textit{New number fields with known \(p\)-class tower}, 
Tatra Mt. Math. Pub.
\textbf{64}
(2015),
21--57, \\
DOI 10.1515/tmmp-2015-0040,
Special Issue on Number Theory and Cryptology \lq 15.
%
\bibitem{Ma2016a}
D. C. Mayer,
\textit{Artin transfer patterns on descendant trees of finite \(p\)-groups}, 
Adv. Pure Math.
\textbf{6}
(2016),
No. 2,
66--104,
DOI 10.4236/apm.2016.62008, 
Special Issue on Group Theory Research,
January 2016.
%
\bibitem{Ma2016b}
D. C. Mayer,
\textit{\(p\)-Capitulation over number fields with \(p\)-class rank two},
J. Appl. Math. Phys.
\textbf{4}
(2016),
No. 7,
1280--1293,
DOI 10.4236/jamp.2016.47135.
%
\bibitem{Ma2018}
D. C. Mayer,
\textit{Modeling rooted in-trees by finite \(p\)-groups},
Chapter 5, pp. 85--113,
in the Open Access Book \textit{Graph Theory --- Advanced Algorithms and Applications},
Ed. B. Sirmacek,
InTech d.o.o., Rijeka, January 2018,
DOI 10.5772/intechopen.68703.
%
\bibitem{Ma2020a}
D. C. Mayer,
\textit{Pattern recognition via Artin transfers applied to class field towers},
3rd International Conference on Mathematics and its Applications (ICMA) 2020,
Facult\'e des Sciences d' Ain Chock Casablanca (FSAC), Universit\'e Hassan II,
Casablanca, Morocco, invited keynote February 28, 2020, \\
\texttt{http://www.algebra.at/DCM@ICMA2020Casablanca.pdf}.
%
\bibitem{Ma2020b}
D. C. Mayer,
\textit{Schur \(\sigma\)-groups with abelian quotient invariants \((9,3)\)},
arXiv:2006.09177.
%
\bibitem{Ma2021a}
D. C. Mayer,
\textit{Bicyclic commutator quotients with one non-elementary component},
arXiv:2108.10754.
%
\bibitem{Ma2021b}
D. C. Mayer,
\textit{BCF-groups with elevated rank distribution},
arXiv:2110.03558.
%
\bibitem{Ma2021c}
D. C. Mayer,
\textit{First excited state with moderate rank distribution},
arXiv:2110.06511.
%
%
\bibitem{Nm1977}
M. F. Newman,
\textit{Determination of groups of prime-power order},
pp. 73--84
in: Group Theory, Canberra, 1975,
Lecture Notes in Math.,
Vol. \textbf{573}
(1977),
Springer,
Berlin.
%
\bibitem{Ob1990}
E. A. O'Brien, 
\textit{The p-group generation algorithm}, 
J. Symbolic Comput.
\textbf{9}
(1990),
677--698.
%
%
\bibitem{Sh1964}
I. R. Shafarevich,
\textit{Extensions with prescribed ramification points} (Russian),
Publ. Math., Inst. Hautes \'Etudes Sci. 
\textbf{18}
(1964),
71--95.
(English transl. by J. W. S. Cassels in
Amer. Math. Soc. Transl.,
II. Ser.,
\textbf{59}
(1966),
128--149.)
%
%
\end{thebibliography}
\end{document}